\documentclass[final,3p,times]{elsarticle}
    \makeatletter
    \def\ps@pprintTitle{%
       \let\@oddhead\@empty
       \let\@evenhead\@empty
       \let\@oddfoot\@empty
       \let\@evenfoot\@oddfoot
    }
    \makeatother

\usepackage{amssymb}
\usepackage[fleqn]{amsmath}
\usepackage{color}
\usepackage{float}

\newtheorem{thm}{Theorem}
\newtheorem{lem}[thm]{Lemma}
\newdefinition{rmk}{Remark}
\newproof{pf}{Proof}
\newdefinition{algo}{Algorithm}
\newdefinition{assumption}{Assumption}
\newcommand{\R}{\mathbb R}
\newcommand{\goto}{\rightarrow}
\newcommand{\f}{\frac}
\newcommand{\nr}[1]{\left\|{#1}\right\|}
\newcommand{\eps}{\varepsilon}
\newcommand{\wit}{ \text{ with }\,}
\DeclareMathOperator{\dd}{\,d}
\DeclareMathOperator{\argmin}{argmin}

\begin{document}

\begin{frontmatter}


\title{Benchmarking results for the Newton-Anderson method}

\author[label1]{Sara Pollock}
\author[label1]{Hunter Schwartz}
\address[label1]
{Department of Mathematics, University of Florida, Gainesville, FL, 32611}

\begin{abstract}
This paper primarily presents numerical results for the Anderson accelerated
Newton method on a 
set of benchmark problems. The results demonstrate superlinear convergence to 
solutions of both degenerate and nondegenerate problems.
The convergence for nondegenerate problems is also justified theoretically.
For degenerate problems, those whose Jacobians are singular at a solution, the 
domain of convergence is studied. It is observed in that setting that 
Newton-Anderson has a domain of convergence similar to Newton, 
but it may be attracted to a different solution than Newton if the problems are slightly
perturbed.
\end{abstract}

\begin{keyword}
Anderson acceleration \sep Newton iteration \sep degenerate problems \sep 
superlinear convergence

\MSC[2010]{65B05,65H10}

\end{keyword}

\end{frontmatter}

\section{Introduction}
\label{sec:intro}
The purpose of this paper is to present some standard benchmark tests, including both 
degenerate and nondegenerate problems, for Anderson acceleration applied to Newton 
iterations, referred to as the Newton-Anderson method. 
The motivation behind studying this method stems from the results of 
\cite{PR19}, where it was shown how Anderson acceleration can locally improve the
convergence rate of linearly converging fixed-point iterations for nondegenerate 
problems, and from the numerical results of \cite{EPRX19}, where Anderson acceleration
applied to Newton's method was demonstrated to allow convergence for a finite 
element discretization of the steady Navier-Stokes equations with a Reynold's number
high enough to cause standard Newton iterations to diverge.  In that setting, the 
convergence history for Newton-Anderson was similar to that of a damped 
Newton method.

A further investigation in \cite{Pollock19} proved that in one dimension, 
Anderson($1$), Anderson acceleration with an algorithmic depth of $m=1$, 
applied to Newton provides higher order (superlinear) convergence to 
nonsimple roots of scalar equations. Newton alone is known to only provide linear 
convergence for such problems \cite[Chapter 6]{Quarteroni06}, 
at the rate $1-1/p$, where $p$ is the multiplicity of the root.  
Generalizing to higher dimensions, it is reasonable to ask whether
Newton-Anderson also provides superlinear convergence for degenerate systems: those
whose Jacobians are singular at a solution.  
The convergence of Newton's method for such systems has been shown to be locally 
linear, in the intersection of a ball and a
star-shaped domain about the solution \cite{Reddien78,Reddien79}.
The numerical results in subsection \ref{subsec:degen} indicate that Newton-Anderson
can provide such locally superlinear convergence, and that an algorithmic depth of
$m=1$ can be sufficient to accomplish this. 

Results for the first accelerated Newton method of \cite{KeSu83},  
which has been shown both theoretically and numerically to 
provide superlinear convergence for degenerate problems \cite{KeSu83,DeKe85,Griewank85},
are shown alongside Newton and Newton-Anderson. 
This method, which can be viewed either 
as a predictor-corrector or an extrapolation method applied to every other iteration,
is similar in form to Anderson(1) (particularly when considered 
from the extrapolation viewpoint), but has convergence properties that are sensitive
to the tuning of two parameters. 
This method will be referred to as accelerated-Newton
(or KS-acc. N., in the tables of results). 
The advantage of Newton-Anderson(1), is that, 
to the present authors' knowledge, there is not a robust method of determining 
{\em a priori}
the two parameters, $C > 0$ ($C \in \R$ as originally presented in \cite{KeSu83})
and $\alpha  \in (0,1)$, of the accelerated-Newton method.
A brief discussion of convergence theory for Newton-Anderson in the nondegenerate
setting is included in subsection \ref{subsec:thy}.
The analysis of the convergence properties of Newton-Anderson however, 
particularly in the degenerate setting, remains pending.

In Section \ref{sec:bench}, the Newton-Anderson method is first tested on problems
from the standard benchmark problem set of \cite{MGH81} 
to identify further problem classes on which Newton-Anderson is 
either advantageous or disadvantageous; and, to identify problems or problem
classes where Newton-Anderson($m$), with $m > 1$, is preferable to Newton-Anderson(1).
The only problem identified from this set that falls in the last category is the Powell 
badly scaled function.  
Another such example designed to demonstrate that there exist problems where depth
$m > 1$ is preferred is \eqref{eqn:and-exact} shown in subsection \ref{subsec:degen}.
That degenerate problem was purposely designed to demonstrate Newton-Anderson($m$)
converging at at iterate $x_{m+1}$, for a problem $f(x) = 0$, of the form 
$f_i(x) = [(Ax - b)_i]^{p_i}$, where $f:\R^n \goto \R^n$, and there are $m \le n$ distinct
exponents $p_i$.

Next, the iterative schemes tested on the problems of sections \ref{sec:bench}
and \ref{sec:dofc} are specified.
Each one is given below in terms of a uniform damping parameter
$\beta \in (0,1]$, as is standard practice for all but Algorithm \ref{alg:KS}.
Damping ($\beta \ne 1$) 
was only used here for one problem, the Brown almost-linear function on $\R^n$,
with $n=20$.
In that case the damping factor $\beta = 0.8$ improved the convergence of all
methods tested. 

\begin{algo}[Newton's method]\label{alg:newton} Set $\beta \in (0,1]$. Choose $x_0$. 
\\ \noindent
For $k = 0, 1, \ldots$
\\ \indent
Set $x_{k+1} = x_k - \beta [f'(x_k)]^{-1}f(x_k)$
\end{algo}

\begin{algo}[Newton-Anderson(1) from \cite{Anderson65}]\label{alg:ander} 
Set $\beta \in (0,1]$. Choose $x_0$. 
Compute $w_1 = -[f'(x_0)]^{-1}f(x_0)$. Set $x_1 = x_0 +\beta w_1$
\\ \noindent
For $k = 1, 2, \ldots$
\\ \indent
Compute $w_{k+1}= -[f'(x_k)]^{-1}f(x_k)$
\\ \indent
Compute $\gamma^{k+1} = {(w_{k+1}, w_{k+1}-w_k)}/{\nr{w_{k+1}- w_k}^2}$
\\ \indent
Set $x_{k+1} = x_k + \beta w_{k+1} - \gamma^{k+1}
\left((x_{k+1}-x_k) + \beta (w_{k+1}-w_k) \right)$
\end{algo}

\begin{algo}[Newton-Anderson($m$) from \cite{Anderson65}]\label{alg:anderm} 
Set $\beta \in (0,1]$ and $m \ge 0$. Choose $x_0$. 
Compute $w_1 = -[f'(x_0)]^{-1}f(x_0)$. Set $x_1 = x_0 + \beta w_1$
\\ \noindent
For $k = 1, 2, \ldots$, set $m_k = \min\{k,m\}$
\\ \indent
Compute $w_{k+1}= -[f'(x_k)]^{-1}f(x_k)$
\\ \indent
Set $F_k= \begin{pmatrix}(w_{k+1}-w_k) & \ldots & (w_{k-m+2} - w_{k-m+1}\end{pmatrix}$, 
and $E_k= \begin{pmatrix}(x_{k}-x_{k-1}) & \ldots & (x_{k-m+1} - x_{k-m})\end{pmatrix}$\\ \indent
Compute $\gamma^{k+1} = \argmin_{\gamma \in \R^m} \nr{w_{k+1} - F_k \gamma}$
\\ \indent
Set $x_{k+1} = x_k + \beta w_{k+1} - \left(E_k + \beta F_k \right)\gamma^{k+1}$
\end{algo}

\begin{algo}[Accelerated Newton from Theorem 1.3 of \cite{KeSu83}]\label{alg:KS}
Set $\beta \in (0,1]$ and parameters $C > 0$ and $\alpha \in (0,1)$. Choose $x_0$.
\\ \noindent
For $k = 1, 2, \ldots$
\\ \indent
Compute $\hat w_{k+1}= -[f'(x_k)]^{-1}f(x_k)$. 
\\ \indent
Set $y_{k+1} = x_k + \beta \hat w_{k+1}$
\\ \indent
Compute $w_{k+1}= -[f'(y_{k+1})]^{-1}f(y_{k+1})$
\\ \indent
Set  $x_{k+1}= y_{k+1} - (2-C\nr{w_{k+1}}^\alpha) w_{k+1}$
\end{algo}
In accordance with \cite{kelley86}, results are shown for Algorithm \ref{alg:KS} 
with parameter choices $C=1$ and $\alpha = 0.9$.
In the cases where the iteration 
with those parameters failed to converge, results are additionally 
shown with parameters 
$C = 0.35$ and $\alpha = 0.1$. In \cite{KeSu83}, where Algorithm \ref{alg:KS} 
is introduced, 
the given parameter range for $C$ is $\R$; however, $C > 0$ appears to be standard in 
presented demonstrations of the method
\cite{KeSu83,DeKe85,Griewank85,kelley86}, and $C \le 0$ is
not considered here.

If the norm $\| \cdot \|$ used in the optimization step of Algorithm \ref{alg:anderm}
is induced by an inner product $(\,\cdot , \cdot \,)$, 
Algorithm \ref{alg:anderm} reduces to Algorithm \ref{alg:ander}, for $m=1$ 
(and reduces to Algorithm \ref{alg:newton}, for $m=0$).  
See \cite{FaSa09,WaNi11} for the equivalence of this
form of the Anderson algorithm to that originally stated in \cite{Anderson65},
and \cite{ToKe15} for a results on optimization in norms not induced by inner
products. Throughout the remainder, the norm $\| \cdot \|$ will denote the $l_2$
norm over $\R^n$.

\subsection{Convergence theory}
\label{subsec:thy}
Convergence and acceleration theory for Anderson acceleration applied to fixed-point
iterations has been recently been studied under assumptions that the underlying
fixed-point operator is contractive (linearly converging) 
\cite{PR19,EPRX19,ToKe15,K18,PRX19}, nonexpansive
\cite{ZOB18}, or nondegenerate \cite{PR19}.  
Here, a small contribution to the theoretical understanding of the method is  
introduced for the particular case of applying Anderson acceleration 
to Newton iterations.
The first lemma shows for nondegenerate problems 
that Newton-Anderson($m$) displays superlinear 
convergence under the assumption  that the optimization 
coefficients are bounded.  
The second lemma shows for Newton-Anderson(1) that if the direction cosine between
consecutive update steps is bounded below unity (the update steps don't become
parallel), then the resulting optimization coefficient is bounded.
This result is applied to the first lemma in Theorem \ref{thm:1}
to show superlinear convergence of the error for Newton-Anderson(1). 
The results that follow in this section do not show that Newton-Anderson is
an improvement over Newton (and where Newton converges quadratically, it generally
is not), but they do show that the technique is reasonable.  

While the underlying theory behind the superlinear convergence in degenerate settings 
is of interest to the authors, it is only demonstrated numerically here.
The first part of the next assumption is a nondegeneracy condition which 
requires that the singular values of $f'$ are uniformly bounded away from zero.
\begin{assumption}\label{assume:f} 
Let $f: \R^n \goto \R^n$, be continuously Frech\'et differentiable,
and suppose there exist positive constants $\sigma$ and $L$ such that
\begin{align}\label{eqn:fnondegen}
\sigma\nr{x} &\le \nr{f'(y)x}, ~\text{ for all }~ x,y \in\R^n, 
\\ \label{eqn:flip}
\nr{f'(x) - f'(y)} &\le L\nr{x-y}, ~\text{ for all}~ x,y \in \R^n.
\end{align}
\end{assumption}
Based on Assumption \ref{assume:f}, $f'$ is invertible on $\R^n$, and since $\sigma$
provides a lower bound on the smallest singular value of $f'$, its inverse 
satisfies
\begin{align*}
\nr{f'(x)^{-1}y} \le \f 1{\sigma} \nr{y}, ~\text{ for all } x,y \in \R^n.
\end{align*}

For the remainder of the section, suppose $f(x) = 0$ has solution $x^\ast$.
The error in iteration $k$ will be denoted by $e_k^\ast = x_k - x^\ast$.

\begin{lem}[Convergence assuming bounded coefficients]\label{lem:1}
Let Assumption \ref{assume:f} hold, and let $x^\ast$ be a solution to the problem
$f(x) = 0$. Then the errors $\{e_{k+1}^\ast\}_{k \ge 1}$
generated by Algorithm \ref{alg:anderm} with $\beta = 1$, satisfy 
\begin{align}\label{eqn:lem1}
\nr{e_{k+1}^\ast} \le \f{L\bar \gamma_{k+1}}{2\sigma} \sum_{j = k-m}^k\nr{e_j^\ast}^2,
\quad
~\text{ where }~
\bar \gamma_{k+1} = \max \left\{2\nr{\gamma^{k+1}}_\infty, 1 + \nr{\gamma^{k+1}}_\infty
\right\}.
\end{align}
For the simple but useful case of depth $m=1$, $\gamma^{k+1}$ is a scalar coefficient
and 
\begin{align}\label{eqn:lem1m1}
\nr{e_{k+1}^\ast} \le \f{L}{2\sigma} \left( |1 - \gamma^{k+1}|\nr{e_k^\ast}^2 
+ |\gamma^{k+1}|\nr{e_{k-1}^\ast}^2\right).
\end{align}
\end{lem}
The results of Lemma \ref{lem:1} show superlinear convergence if 
$\bar \gamma_{k+1}$ (or $\gamma^{k+1}$ for $m=1$) can be shown to be bounded.
In Lemma \ref{lem:2}, 
a condition is found on the direction cosines between $w_{k+1}$ and $w_k$ to 
guarantee this bound for the particular case of depth $m=1$.
While the more general result
for depths $m > 1$ is not presented here, it can likely be obtained through an
application of the 
machinery developed in \cite{PR19} in which a QR decomposition is used to analyze 
the defect introduced by the matrix $F_k$ of Algorithm \ref{alg:anderm} 
having columns that are not orthogonal.
\begin{pf}
For an arbitrary index $l \ge 1$, the usual estimates for Newton's method hold 
under Assumption \ref{assume:f}. In particular expanding $f(x_l)$ about 
$0 = f(x^\ast)$ for an arbitrary index $l \ge 0$ allows
$f(x_l)  = \int_0^1 f'(x^\ast + t e_l^\ast) e_l^\ast \dd t$.  
The error in $x_l + w_{l+1}$ then satisfies
\begin{align*}
x_l + w_{l+1} - x^\ast & = e_l^\ast - f'(x_l)^{-1} f(x_l)
= f'(x_l)^{-1}(f'(x_l) e_l^\ast - f(x_l))
= f'(x_l)^{-1} \int_0^1 (f'(x_l)  - f'(x^\ast + t e_l^\ast)) e_l^\ast \dd t.
\end{align*}
Under the two conditions of Assumption \ref{assume:f}, the usual bound holds
\begin{align}\label{eqn:l1-002n}
\nr{x_l + w_{l+1} - x^\ast} \le \f{L}{2 \sigma} \nr {e_l^\ast}^2.
\end{align} 

Expanding the update step of Algorithm \ref{alg:anderm}, with 
$m_k$ denoted by $m$, and the entries of $\gamma^{k+1}$ denoted 
$\gamma_1, \ldots, \gamma_m$, for simplicity
\begin{align*}
x_{k+1} &= x_k + w_{k+1} - \sum_{j = k-m+1}^k((x_j-x_{j-1}) + (w_{j+1} - w_j))
\gamma_{k-j+1}
\nonumber \\
& = (x_k + w_{k+1})(1-\gamma_1) + (x_{k-1}+ w_k)(\gamma_1 - \gamma_2)
+ \ldots + (x_{k-m+1} + w_{k-m+2})(\gamma_{m-1} - \gamma_m) 
+ (x_{k-m} + w_{k-m+1})\gamma_m
\nonumber \\
& = (x_k + w_{k+1})(1-\gamma_1) + \sum_{j = k-m+1}^{k-1}(x_{j}+ w_{j+1})
(\gamma_{k-j} - \gamma_{k-j+1})
+ (x_{k-m} + w_{k-m+1})\gamma_m.
\end{align*}
Writing $e^\ast = e^\ast((1 - \gamma_1) + (\gamma_1 - \gamma_2) + \ldots
+ (\gamma_{m-1} - \gamma_m) + \gamma_m)$, the error $e_{k+1}^\ast$ is represented by
\begin{align}\label{eqn:l1-004}
e^\ast_{k+1} = (x_k + w_{k+1} - x^\ast)(1-\gamma_1) 
+ \sum_{j = k-m+1}^{k-1}(x_{j}+ w_{j+1} - x^\ast)
(\gamma_{k-j} - \gamma_{k-j+1})
+ (x_{k-m} + w_{k-m+1} - x^\ast)\gamma_m.
\end{align}
Taking norms and applying \eqref{eqn:l1-002n} to each term of \eqref{eqn:l1-004},
yields the result \eqref{eqn:lem1}.

For the simple case of $m=1$, letting $\gamma = \gamma^{k+1}$, the estimate
\eqref{eqn:l1-002n} is applied to the expansion
\[
x_{k+1} - x^\ast = x_k + w_{k+1} - x^\ast - 
\gamma\left((x_k - x_{k-1}) + (w_{k+1} - w_k) \right)
= (1-\gamma)(x_k + w_{k+1} - x^\ast) + \gamma(x_{k-1} + w_k - x^\ast),
\]
to produce the result \eqref{eqn:lem1m1}.
\end{pf} 

The next lemma is specific to the algorithmic depth of $m=1$, and produces a simple
verifiable condition under which the optimization coefficient $\gamma^{k+1}$ for
depth $m=1$ is bounded in magnitude by 1.
\begin{lem}[Boundedness conditions for $\gamma^{k+1}$ for Newton-Anderson(1)]
\label{lem:2}
The coefficient $\gamma^{k+1}$ produced by Algorithm \ref{alg:ander}
(or Algorithm \ref{alg:anderm} with $m=1$)
is bounded, so long as $\cos(w_{k+1},w_k) \le \bar \alpha$, for some $\bar \alpha < 1$.
Specifically, for $\bar \alpha = 0.942$, it holds that $|\gamma| < 2$.
\end{lem}
\begin{pf}
For depth $m=1$, $\gamma = \gamma^{k+1}$ is given by 
\begin{align}\label{eqn:lem2-001}
\gamma = \f{(w_{k+1},w_{k+1}-w_k)}{\nr{w_{k+1}- w_k}^2}
= \f{\nr{w_{k+1}^2} - (w_{k+1},w_k)}{ \nr{w_{k+1}}^2 + \nr{w_{k}}^2 - 2(w_{k+1},w_k)}
= \f{\nr{w_{k+1}}/\nr{w_k} - \cos(w_{k+1},w_k)}
    {\nr w_{k+1}/\nr{w_k} + \nr{w_k}/\nr{w_{k+1}} - 2\cos(w_{k+1},w_k)},
\end{align}
where $\cos(u,v)$ is the direction cosine between vectors $u,v \in\R^n$.
For $w_{k+1}, w_k \ne 0$, the right-hand side of \eqref{eqn:lem2-001} is of the form 
\begin{align}\label{eqn:lem2-002}
\gamma = \psi_\alpha(r)
= \f{r - \alpha}{r + 1/r - 2\alpha}, \quad r >0, ~-1 \le \alpha \le 1. 
\end{align}

Clearly $\psi_{1}(r)$ has a singularity as $r \goto 1$; however for $\alpha < 1$,
there are no singularities.  Hence, so long as
$\cos(w_{k+1},w_k) \le \bar \alpha < 1$, the coefficient $\gamma$ remains bounded.  

To be more precise and determine a value of $\bar \alpha$ for which $|\gamma| < 2$,
for instance, 
it is first noted that $|\psi_\alpha(r)| < 1$ for $\alpha \in [-1,0]$ and 
$r>0$, and $\psi_\alpha(1) = 1/2$.
So, it remains to investigate $\alpha \in (0,1)$ for $r > 0, ~r \ne 1$.
The extreme values of $\psi_\alpha(r)$ can first be found for each fixed value of 
$\alpha \in (0,1)$, by setting 
$0 = \psi_\alpha'(r) = (-\alpha r^2 + 2r - \alpha)/(r^2 - 2\alpha r+1)^2$.
The numerator yields a quadratic equation in $r$ 
to which the solution for $r < 1$ is 
$r_\alpha = (1 - \sqrt{1-\alpha^2})/\alpha$. 
Plugging this into \eqref{eqn:lem2-002}
yields
\begin{align*}
\psi_\alpha(r_\alpha) = \f{(2/\alpha^2-1)(1 - \sqrt{1 - \alpha^2})-1}
{2(1/\alpha^2 -1)(1 - \sqrt{1 - \alpha^2})},
\end{align*}
which is a decreasing function of $\alpha$, for $\alpha > 0$. Investigating this
function numerically, it is seen, for example, that $\psi_\alpha(r) > -0.990$ for 
$\alpha < 0.942$.  

For $r > 1$, the extremum of $\psi_\alpha(r)$ occurs at
$r_\alpha = (1 + \sqrt{1-\alpha^2})/\alpha$, for which
\begin{align*}
\psi_\alpha(r_\alpha) = \f{(2/\alpha^2-1)(1 + \sqrt{1 - \alpha^2}) -1}
{2(1/\alpha^2 -1)(1 + \sqrt{1 - \alpha^2})},
\end{align*}
which is an increasing function of $\alpha$, for $\alpha > 0$; and, for which 
$\psi_\alpha(r) < 1.99$ for $\alpha < 0.942$.

\end{pf}

\begin{thm}[Superlinear convergence for Newton-Anderson(1)]\label{thm:1}
Let Assumption \ref{assume:f} hold, and let $x^\ast$ be a solution to the problem
$f(x) = 0$. Then the errors $\{e_{k+1}^\ast\}_{k \ge 1}$
generated by Algorithm \ref{alg:ander} (or Algorithm \ref{alg:anderm} with $m=1$) 
and $\beta = 1$
converge superlinearly to zero if the direction cosines between update steps satisfy
$\cos(w_{k+1},w_k) < \bar \alpha$, for some constant $\bar \alpha < 1$.
\end{thm}
\begin{pf}
The result follows immediately from applying the results of Lemma \ref{lem:2} to 
the result \eqref{eqn:lem1m1} of Lemma \ref{lem:1}.
\end{pf}
This shows that if consecutive update steps are not parallel 
(anti-parallel is harmless), then convergence is locally better than linear.
The superlinear convergence for nondegenerate problems agrees with the numerical 
results shown below. Consistent with the theory, the convergence is  
also observed to be subquadratic.  
While it is not demonstrated here, superlinear convergence
for degnererate problems (which violate Assumption \ref{assume:f}) is 
also obesrved.

\section{Benchmark results}
\label{sec:bench}
The first set of problems considered are the systems of nonlinear equations
that map $f: \R^n \goto \R^n$, from the standard test set of Mor\'e {\em et al.},
\cite{MGH81}. The 
problems and results are stated in subsection \ref{subsec:MGH}.
Subsection \ref{subsec:degen} 
that follows specifically considers additional degenerate problems.

Throughout this section, each method was terminated on the residual reduced
below the tolerance $\|f(x_k)\| < 10^{-8}$.  In section \ref{sec:dofc} discussing
the domain of convergence to particular solutions, iterations
were terminated based on error dropping below tolerance $\|x_k - x^\ast\| < 10^{-8}$,
where $x^\ast$ is a given (known) solution.
All computations were done in Matlab on an
8 core Intel Xeon W-2145 CPU @ 3.70GHz, with 64 GB memory.
Timing was performed using Matlab's tic/toc commands. 
Each algorithm/problem pair was run 100 times, 
and the results were averaged to produce the time displayed in the tables of results. 

Problems from the {\em Systems of Nonlinear Equations} test set of \cite{MGH81} for
which all methods converged in no more than four iterations, in addition to  the
Chebyquad function, are excluded from the test set below.
In a slight (but standard) abuse of notation, subscripts are used in the algorithms 
above to indicate iteration counts, and in the problem descriptions below to 
denote components of $x \in \R^n$.

\subsection{Problems $f: \R^n \goto \R^n$ from \cite{MGH81}}\label{subsec:MGH}
\begin{enumerate}
\item[B1.] Powell badly scaled function, \cite[(3)]{MGH81} and \cite{Powell70}, n = 2.
\begin{align}\label{eqn:pow-bd-sc}
f(x) = \begin{pmatrix}
10^4 x_1 x_2 -1 \\ \exp(-x_1) + \exp(-x_2) -1.0001
\end{pmatrix}, \quad
\wit x_0 = \begin{pmatrix} 0 , 1\end{pmatrix}^T.
\end{align}

\item[B2.] Helical valley function, \cite[(7)]{MGH81} and \cite{FlPo63}, n = 3.
\begin{align}\label{eqn:helicalv}
f(x) = \begin{pmatrix}
10(x_3 - 10 \theta(x_1, x_2))\\
10( \sqrt{x_1^2 + x_2^2} -1) \\
x_3
\end{pmatrix}, 
\quad \text{ where } \theta(x_1, x_2) = \f 1 {2 \pi}\arctan(x_2/x_1),
\quad
\wit x_0 = \begin{pmatrix} -1 , 0 , 0\end{pmatrix}^T,
\end{align}
where $\arctan(\cdot)$ is the four-quadrant arctangent function.

\item[B3.] Powell singular function, \cite[(13)]{MGH81} and \cite{powell62}, n = 4. 
\begin{flalign}\label{eqn:pow-sing}
f(x) &= \begin{pmatrix} 
x_1 + 10 x_2 \\ \sqrt 5 (x_3 - x_4) \\ (x_2 - 2x_3)^2 \\ \sqrt{10}(x_1 - x_4)^2
\end{pmatrix}, 
\wit x_0 = \begin{pmatrix} 3 , -1 , 0 , 1 \end{pmatrix}^T.
\quad
\end{flalign}

\item[B4.] Watson function, \cite[(20)]{MGH81} and \cite{KoOs68}, $n = 31$.
\begin{flalign}\label{eqn:watson}
f_i(x) &= \sum_{j = 2}^n (j-1)x_j t_i^{j-2} 
 - \left(\sum_{j=1}^n x_j t_i^{j-1}\right)^2 -1, 
\text{ where }t_i = i/29, ~1 \le i \le 29, 
\nonumber \\
f_{30}(x) &= x_1, ~f_{31}(x) = x_2 - x_1^2-1,
\quad
\wit x_0 = \begin{pmatrix}0, \ldots, 0 \end{pmatrix}^T.
\end{flalign}

\item[B5.] Trigonometric function, \cite[(26)]{MGH81} and \cite{FlPo63}, 
$n = 100,1000,10000$. 
\begin{align}\label{eqn:trig}
f_i(x) = n - \sum_{j = 1}^n \cos x_j + i(1 - \cos x_i ) - \sin x_i,
\quad
\wit x_0 = \begin{pmatrix}1/n, \ldots, 1/n \end{pmatrix}^T.
\end{align}

\item[B6.] Brown almost-linear function, \cite[(27)]{MGH81} and \cite{Brown69}, 
$n = 5,20$.
\begin{flalign}\label{eqn:brown-al-li}
f_i(x) &= x_i + \sum_{j = 1}^n x_j - (n+1), \quad 1 \le i < n, 
\nonumber \\
f_n(x) &= \left(\prod_{j=1}^n x_j\right) -1,
\quad
\wit x_0 = \begin{pmatrix}1/2, \ldots, 1/2 \end{pmatrix}^T.
\end{flalign}

\item[B7.] Broyden tridiagonal function, \cite[(30)]{MGH81} and \cite{broyden65}, 
$n = 1000$.
\begin{align}\label{eqn:broyden-tridiag}
f_i(x) = (3-2x_i)x_i - x_{i-1}- 2x_{i+1} + 1, ~(x_0 = x_{n+1}= 0),
\quad
\wit x_0 = \begin{pmatrix}-1, \ldots, -1 \end{pmatrix}^T.
\end{align}

\item[B8.] Broyden banded function, \cite[(31)]{MGH81} and \cite{broyden71}, 
$n = 1000$.
\begin{align}\label{eqn:broyden-bdd}
f_i(x) &= x_i(2 + 5 x_i^2) + 1 - \sum_{j \in J_i}x_j(1+x_j),
~J_i = \{j\,: j\ne i, \, \max(1, i-5) \le j \le \min(n, i+ 1)\},
\nonumber \\
&\wit x_0 = \begin{pmatrix}-1,\ldots,-1 \end{pmatrix}^T.
\end{align}

\end{enumerate}

\subsubsection{Results}
The first set of results, those for which Anderson(1) was preferred over Anderson($m$), 
with $m > 1$, are summarized in table \ref{tab:A1con}.  Here, results are given
for Newton, Newton-Anderson(1), and accelerated-Newton. An entire predictor-corrector 
step is counted as a single iteration for the accelerated-Newton method.
For each problem and method, the number of iterations ($k$), terminal
residual $\nr{f(x_k)}$, norm of the final update step $\nr{w_k}$, the 
terminal approximate order of convergence 
$q_k = \log \nr{f(x_k)}/\log \nr{f(x_{k-1})}$, 
and the (average) time in seconds, are displayed.
\begin{table}[ht]
\centering
{\renewcommand{\arraystretch}{1.1}
\begin{tabular}{r||c|c|c|c|c|c}
\hline
Problem & method &iterations ($k$)  & $\|f(x_k)\|$ & $\|w_{k}\|$ &
$q_{k}$ &  time (sec)  \\ [2pt]
\hline
B2. \eqref{eqn:helicalv}
&Newton                  & 10 & 1.325e-14 & 4.161e-08 & 2.235 & 9.655e-05 \\
$n=3$&N. Anderson(1)          & 10 & 5.485e-13 & 1.001e-08 & 1.794 & 0.0001249 \\
&KS-acc. N. $(1.0, 0.9)$ & F & -- & -- & -- & -- \\
&KS-acc. N. $(0.35, 0.1)$ & 11 & 1.042e-16 & 6.591e-18 & 2.391 & 0.0001292 \\

\hline
B3. \eqref{eqn:pow-sing}
&Newton                  & 16 & 2.954e-09 & 3.743e-05 & 1.076 & 0.0002417 \\
$n=4$&N. Anderson(1)          & 3 & 4.163e-17 & 0.2157 & 40.27 & 0.0001412 \\
&KS-acc. N. $(1.0, 0.9)$ & 3 & 5.964e-10 & 0.003071 & 2.638 & 0.0001279 \\

\hline
B4. \eqref{eqn:watson}
&Newton                  & 5 & 6.004e-14 & 13.5 & 1.71 & 0.01759 \\
$n=31$&N. Anderson(1)          & 7 & 2.863e-10 & 3.341 & 1.73 & 0.02455 \\
&KS-acc. N. $(1.0, 0.9)$ & 4 & 2.215e-13 & 9.818 & 2.807 & 0.02716 \\

\hline
B5. \eqref{eqn:trig}
&Newton                  & 10 & 1.892e-11 & 6.137e-07 & 1.726 & 0.003895 \\
$n=100$&N. Anderson(1)          & 8 & 9.565e-13 & 1.159e-08 & 1.518 & 0.003095 \\
&KS-acc. N. $(1.0, 0.9)$ & 7 & 3.568e-12 & 3.571e-12 & 1.741 & 0.003868 \\

\hline
B5. \eqref{eqn:trig} 
&Newton                  & 13 & 9.906e-11 & 4.454e-07 & 1.575 & 0.4408 \\
$n=1000$&N. Anderson(1)          & 11 & 1.653e-11 & 2e-08 & 1.4 & 0.3734 \\
&KS-acc. N. $(1.0, 0.9)$ & 6 & 1.302e-09 & 1.308e-09 & 1.534 & 0.3695 \\

\hline
B7. \eqref{eqn:broyden-tridiag} 
&Newton                  & 4 & 1.065e-09 & 4.555e-05 & 2.312 & 0.01596 \\
$n=1000$&N. Anderson(1)          & 6 & 1.59e-14 & 6.612e-09 & 1.788 & 0.02284 \\
&KS-acc. N. $(1.0, 0.9)$ & 4 & 1.904e-13 & 6.395e-14 & 2.187 & 0.02831 \\

\hline
B8. \eqref{eqn:broyden-bdd} 
&Newton                  & 8 & 5.669e-09 & 2.419e-08 & 1.24 & 0.09291 \\
$n=1000$&N. Anderson(1)          & 9 & 3.668e-10 & 1.768e-09 & 1.2 & 0.1056 \\
&KS-acc. N. $(1.0, 0.9)$ & 9 & 3.204e-10 & 6.07e-11 & 1.19 & 0.1814 \\

\hline 
B6. \eqref{eqn:brown-al-li} 
&Newton                  & 18 & 3.14e-15 & 6.481e-08 & 2.058 & 0.0002524 \\
$n=5$ &N. Anderson(1)          & 24 & 5.031e-12 & 2.935e-07 & 1.555 & 0.0003786 \\
&KS-acc. N. $(1.0, 0.9)$ & F & -- & -- & -- & -- \\
&KS-acc. N. $(0.35, 0.1)$ & 9 & 3.397e-15 & 2.334e-16 & 1.963 & 0.0002039 \\

\hline
B6. \eqref{eqn:brown-al-li}
&Newton, $\beta = 0.8$ & 368 & 4.743e-09 & 4.854e-07 & 1.092 & 0.02743 \\
$n=20$&N. Anderson(1), $\beta = 0.8$          & 52 & 7.394e-10 & 3.853e-09 & 1.088 & 0.004206 \\
&KS-acc. N. $(1.0, 0.9)$, $\beta = 0.8$ & F & -- & -- & -- & -- \\
&KS-acc. N. $(0.35, 0.1)$, $\beta = 0.8$ & 147 & 1.931e-09 & 4.224e-08 & 1.602 & 0.02054 \\

\end{tabular}
\caption{Results for problems from \cite{MGH81} where Newton-Anderson(1) converges }
\label{tab:A1con}
}
\end{table}

\begin{table}[ht]
\centering
{\renewcommand{\arraystretch}{1.1}
\begin{tabular}{r||c|c|c|c|c|c}
\hline
Problem & method & iterations ($k$)  & $\|f(x_k)\|$ & $\|w_{k}\|$ &
$q_k$ &  time (sec)  \\ [2pt]
\hline
B1. \eqref{eqn:pow-bd-sc}
&Newton                  & 12 & 1.573e-11 & 3.987e-05 & 1.769 & 0.0001241 \\
$n=2$&N. Anderson(1)          & F & -- & -- & -- & -- \\
&N. Anderson(2)          & 12 & 4.058e-09 & 0.0001451 & 1.518 & 0.0003492 \\
&KS-acc. N. $(1.0, 0.9)$ & 8 & 5.773e-15 & 5.886e-16 & 1.818 & 7.731e-05 \\

\hline
\end{tabular}
\caption{Results for problems from \cite{MGH81} where Newton-Anderson(1) fails to
  converge.}
\label{tab:A1fail}
}
\end{table}

\subsubsection{Discussion}
\begin{figure}[ht!]
\centering
\includegraphics[trim = 100pt 240pt 100pt 250pt,clip = true, width=0.44\textwidth]
{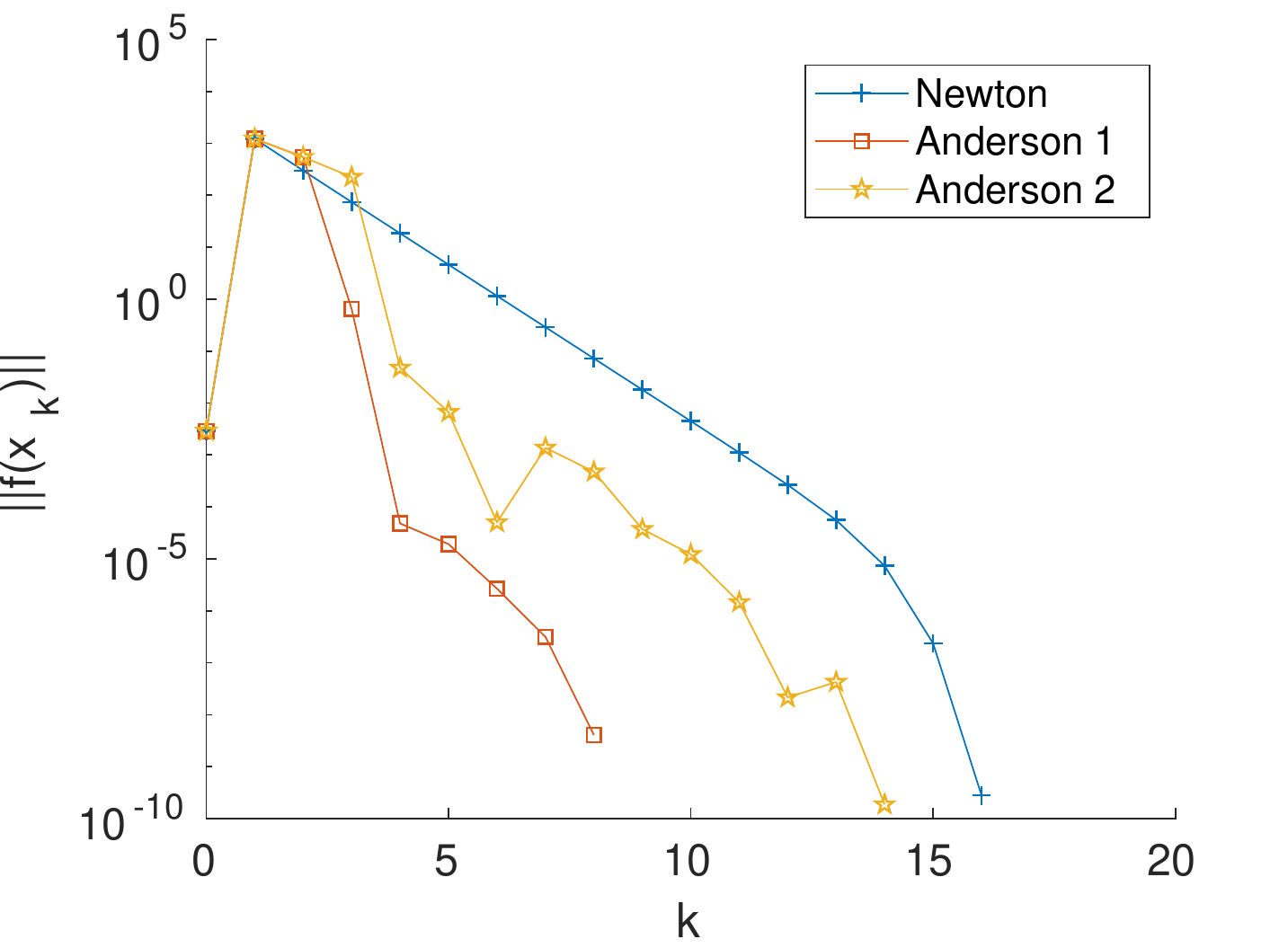}
\includegraphics[trim = 100pt 240pt 100pt 250pt,clip = true, width=0.44\textwidth]
{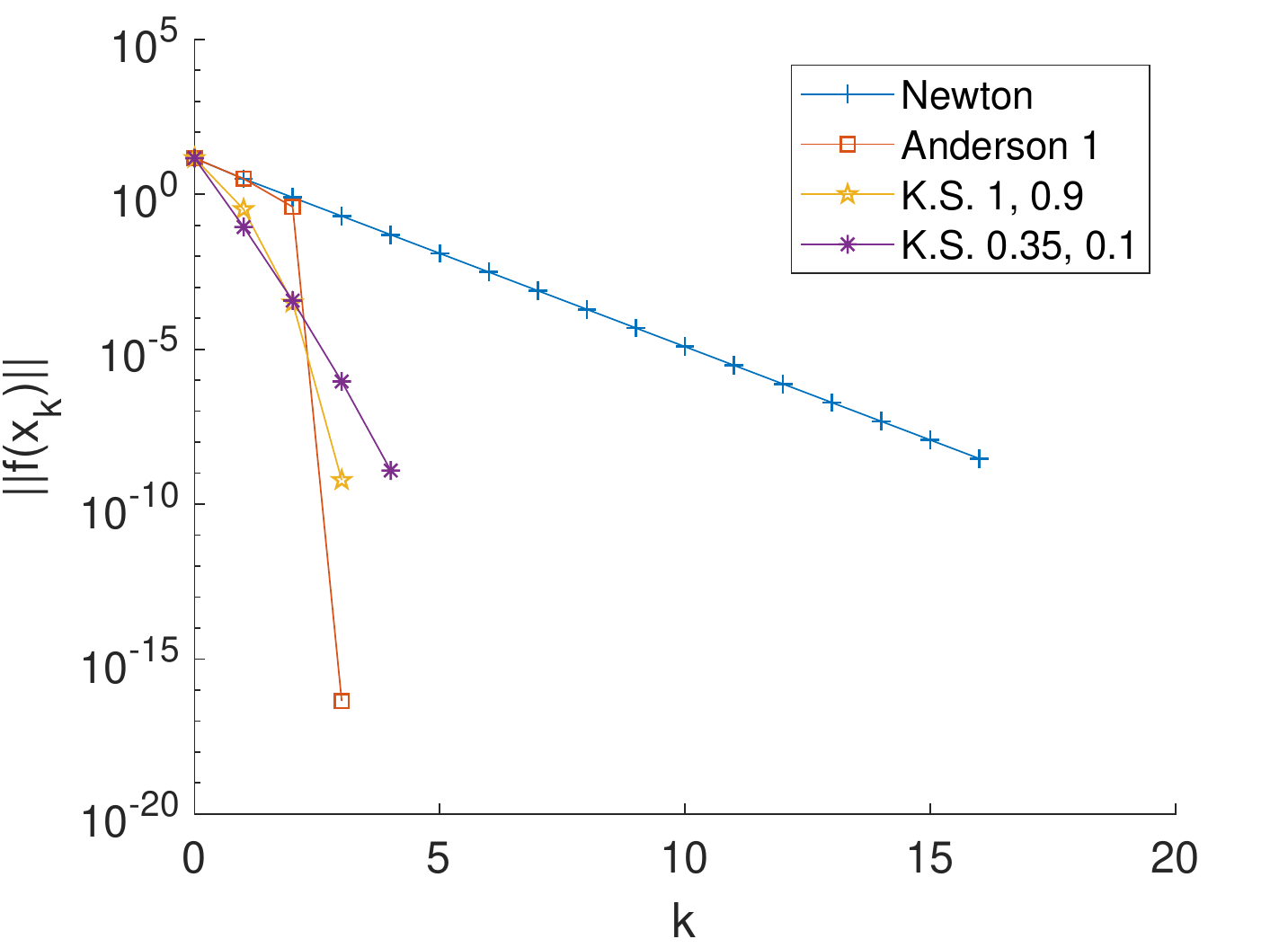}
\caption{Left: convergence history for the trigonometric function \eqref{eqn:trig}, 
with $n=$10,000, for Newton and Newton-Anderson($m$), $m = 1,2$.
Right: convergence hisory for the Powell singular function \eqref{eqn:pow-sing},
illustrating the high terminal order of convergence for Newton-Anderson(1).
\label{fig:trig}
}
\end{figure}

In accordance with the expected results, 
for the Powell singular function \eqref{eqn:pow-sing}, whose Jacobian
is singular at the solution (with rank 2), Newton converges linearly, 
and both accelerated and Anderson methods yield faster convergence, as seen in the 
convergence histories shown in
in figure \ref{fig:trig}, on the right.
In terms of timing, Newton-Anderson(1) is advantageous over accelerated-Newton as
the problem dimension $n$ increases, as Newton-Anderson(1) solves a linear system 
plus computes two additional inner products to set the optimization parameter, whereas
accelerated-Newton solves two linear systems.  This difference is illustrated by 
comparing the results of problem \eqref{eqn:pow-sing} with $n=4$ to problem 
\eqref{eqn:broyden-bdd} with $n=1000$, where both methods converge in the same number
of iterations.

Newton-Anderson(1) was found disdvantageous compared to Newton in terms of timing
for the Watson function \eqref{eqn:watson}, with $n=31$, and the  
Broyden tridiagonal function with $n=1000$ 
\eqref{eqn:broyden-tridiag}, taking respectively 1.40 and 1.44 times as long as Newton 
to converge (but only two additional iterations). It was found nominally 
disadvantageous for the Broyden banded function \eqref{eqn:broyden-bdd}, with $n=1000$,
taking 1.14 times as long to converge as Newton, with one additional iteration.
Newton-Anderson and Newton both converged in 10 iterations for the Helical valley
function \eqref{eqn:helicalv}, with $n=3$.  In each of the above cases, Newton-Anderson
converged faster than accelerated-Newton.  In agreement with the results of
\cite{PR19,EPRX19}, these results indicate that 
Anderson slows the convergence of quadratically-converging
iterations, and should not be used on problems
for which Newton is known to converge robustly and quadratically.

For the Brown almost-linear function \eqref{eqn:brown-al-li}, Newton-Anderson converged
in 1.5 times as long as Newton (6 additional iterations) with $n=5$ and no damping,
but 0.15 times as long as Newton (316 fewer iterations), for the problem with
dimension $n=20$ with a damping factor of $0.8$.  
In both cases, accelerated-Newton failed to converge with parameters 
$(C,\alpha) = (1.0,0.9)$, and succeeded with parameters $(0.35,0.1)$.
Taken together, these results  indicate that Newton-Anderson can be useful for problems 
where convergence of Newton is locally quadratic but less stable as the dimensions of 
the problem increases.

For the trigonometric function \eqref{eqn:trig},
as $n$ increases, the Newton method has an extended preasymptotic regime of linear
convergence. Both accelerated methods, which display superlinear
(but subquadratic) convergence, attain residual tolerance in fewer iterations.
Newton-Anderson($m$) for $m=1,2$ 
are illustrated for $n =$ 10,000, in figure \ref{fig:trig} on the left.
On the right of figure \ref{fig:trig}, the residual history for the Powell singular
function is shown.  The plot illustrates the curiously high approximate order of
convergence for Anderson(1). The accelerated-Newton method is shown for comparison,
with both sets of tested parameters, along with Newton.  

Results for the Powell badly scaled function \eqref{eqn:pow-bd-sc}
are shown in table \ref{tab:A1fail}.
For this problem, Newton-Anderson($m$) converges for $m >1$, but not $m=1$. 
Tests (not shown) were run for $m = 1, \ldots, 20$.  However, as $n=2$ for this
problem, the optimization problem reduces to a solvable linear system for depth
$m=2$, and results with greater depths are not shown.
Convergence can be restored for $m=1$ by applying a safeguarding strategy based on 
Lemma \ref{lem:2}. If a Newton step is taken whenever $\cos(w_{k+1},w_k) > 0.942$,
and the Anderson(1) step otherwise, then Newton-Anderson(1) converges, essentially 
tracking the Newton iteration.

\subsection{Additional degenerate problems}\label{subsec:degen}
Based on the results of \cite{Pollock19}, which show superlinear convergence of
Newton-Anderson(1) to nonsimple roots of scalar equations, it is conjectured 
that Newton-Anderson also provides superlinear convergence to solutions of degenerate 
problems.  Three such problems (in addition to the Powell singular function above)
are collected in this section.  The first two are standard problems from the
literature and the third was chosen to demonstrate a case where Newton-Anderson(4) 
converges faster than Newton-Anderson($m$), for $m = 1,2,3$.

\begin{enumerate}
\item[D1.] From \cite{Reddien79}, $f: \R^3 \goto \R^3$
\begin{align}\label{eqn:poly3indef}
f(x_1,x_2,x_3) = \begin{pmatrix}x_1 + x_1x_2 + x_2^2 \\ x_1^2 - 2x_1 + x_2^2 
                 \\ x_1 + x_3^2
\end{pmatrix}
\quad
\wit x_0 = \begin{pmatrix}0.1,0.5, 1.0 \end{pmatrix}^T.
\end{align}
\item[D2.] Chandrasekhar H-equation from \cite{K18}.
The Chandrasekhar H-equation from radiative transfer
(see \cite{K18,busbridge60,chandra60} and the references therein), is given by
\[
F(H)(\mu) \coloneqq H(\mu) - \left(1 - \f{\omega}{2}\int_0^1 
\f{\mu H(\nu) d\nu}{\mu + \nu} \right)^{-1} = 0.
\]
Discretizing the equation (as described in \cite[Example 2.10]{K18}) by the
composite midpoint rule with $n$ subintervals,
yields the discrete system in $f: \R^n\goto \R^n$ given componentwise for
$h\in \R^n = (h_1, \cdots, h_n)^T$, by
\begin{align}\label{eqn:disc-heqn}
f_i(h) = h_i - \left(1 - \f{\omega}{2n}\sum_{j = 1}^n \f{i h_j}{i + j -1} 
\right), \quad 1 \le i \le n,
\quad
\wit x_0 = \begin{pmatrix}1, \ldots, 1 \end{pmatrix}^T.
\end{align}

\item[D3.] Problem designed to demonstrate Newton-Anderson(4).
\begin{flalign}\label{eqn:and-exact}
f_i(x) &= [(Ax - b)_i]^{p_i}, 
~ A = \begin{pmatrix}2 & -1 &  \\ -1 & 2 & -1 \\
  & \ddots & \ddots & \ddots \\
 & & -1 & 2 & -1\\ & & & -1 & 2
\end{pmatrix}, 
~b = \begin{pmatrix} -11, -7,-5, -3, -2, 2, 3, 5, 7, 11 \end{pmatrix}^T,
\nonumber \\ 
&p = \begin{pmatrix}2, 4, 4, 2, 2, 8, 8, 2, 12, 12 \end{pmatrix}, 
\quad
\wit x_0 = \begin{pmatrix}0, \ldots, 0 \end{pmatrix}^T.
\end{flalign} 
\end{enumerate}

\begin{table}[ht]
\centering
{\renewcommand{\arraystretch}{1.1}
\begin{tabular}{r||c|c|c|c|c|c}
\hline
Problem & method & iterations ($k$)  & $\|f(x_k)\|$ & $\|w_{k}\|$ &
$q_k$ &  time (sec)  \\ [2pt]
\hline
D1.
&Newton                  & 14 & 3.991e-09 & 6.903e-05 & 1.077 & 0.000129 \\
$n=3$ &N. Anderson(1)          & 5 & 1.656e-10 & 1.349e-05 & 1.493 & 7.152e-05 \\
&KS-acc. N. $(1.0, 0.9)$ & 3 & 2.396e-09 & 0.0001877 & 1.993 & 4.476e-05 \\

\hline
D2.
&Newton                  & 16 & 2.628e-09 & 0.000382 & 1.075 & 0.9698 \\
$n=1000$&N. Anderson(1)          & 6 & 1.236e-11 & 0.001947 & 1.663 & 0.3671 \\
&N. Anderson(2)          & 6 & 1.912e-11 & 2.595e-06 & 1.399 & 0.3658 \\
&KS-acc. N. $(1.0, 0.9)$ & F & -- & -- & -- & -- \\
&KS-acc. N. $(0.35, 0.1)$ & 4 & 3.986e-09 & 0.002449 & 1.461 & 0.4772 \\

\hline
D3.
&Newton                  & 46 & 4.339e-09 & 0.07587 & 1.057 & 0.001038 \\
$n=10$&N. Anderson(1)          & 17 & 7.899e-09 & 0.01347 & 1.45 & 0.0008288 \\
&N. Anderson(2)          & 26 & 6.781e-11 & 0.0003507 & 1.404 & 0.002012 \\
&N. Anderson(3)          & 6 & 3.964e-10 & 0.09431 & 5.87 & 0.000219 \\
&N. Anderson(4)          & 5 & 7.289e-25 & 0.1056 & 23.05 & 0.0001861 \\
&KS-acc. N. $(1.0, 0.9)$ & F & -- & -- & -- & -- \\
&KS-acc. N. $(0.35, 0.1)$ & F & -- & -- & -- & -- \\
&KS-acc. N. $(0.7, 0.3)$ & 20 & 1.758e-09 & 0.07474 & 1.165 & 0.0007143 \\

\end{tabular}
\caption{Results for degenerate problems.}
\label{tab:Q6.11-2}
}
\end{table}

\subsection{Discussion}
In problems \eqref{eqn:poly3indef}, \eqref{eqn:disc-heqn}, and 
\eqref{eqn:and-exact},
Newton-Anderson(1) converges faster than Newton by respective factors of
0.55, 0.38, and 0.8. Newton-Anderson(4) outperforms Newton-Anderson(1) by a factor
of 0.22, and accelerated-Newton by a factor or 0.26, in the last problem.
The approximate convergence order $q_k$ is close to one for Newton in each case,
and greater than one (generally above 1.4) for each of the accelerated methods, when
they converge. 
The accelerated-Newton method was run with an additional parameter pair for 
\eqref{eqn:and-exact} as the method did not converge under either of the other 
two parameter-pairs used in the rest of the tests.
The varied nature of these problems in terms of both problem dimension, and
the dimension of the nullspaces at their respective solutions, indicates consistent
behavior of Newton-Anderson on degenerate problems, which is worth a theoretical 
investigation.

The first problem \eqref{eqn:poly3indef} has a Jacobian featuring a two-dimensional
nullspace at the solution.
For the Chandrasekhar $H$-equation \eqref{eqn:disc-heqn}, 
the rank of the Jacobian at the solution is $n-1$.
For the degenerate polynomial problem \eqref{eqn:and-exact}, 
the Jacobian has rank zero at the solution.
The problem was designed to solve exactly for $m = 4$ after the first full 
optimization.  It is shown here to illustrate that there exists a class of 
problems for which $m > 1$ is preferable.  It is not yet clear what this 
problem class is.

\section{Domain of convergence}\label{sec:dofc}
\begin{figure}[ht!]
\centering
\includegraphics[trim = 90pt 240pt 110pt 250pt,clip = true, width=0.32\textwidth]
{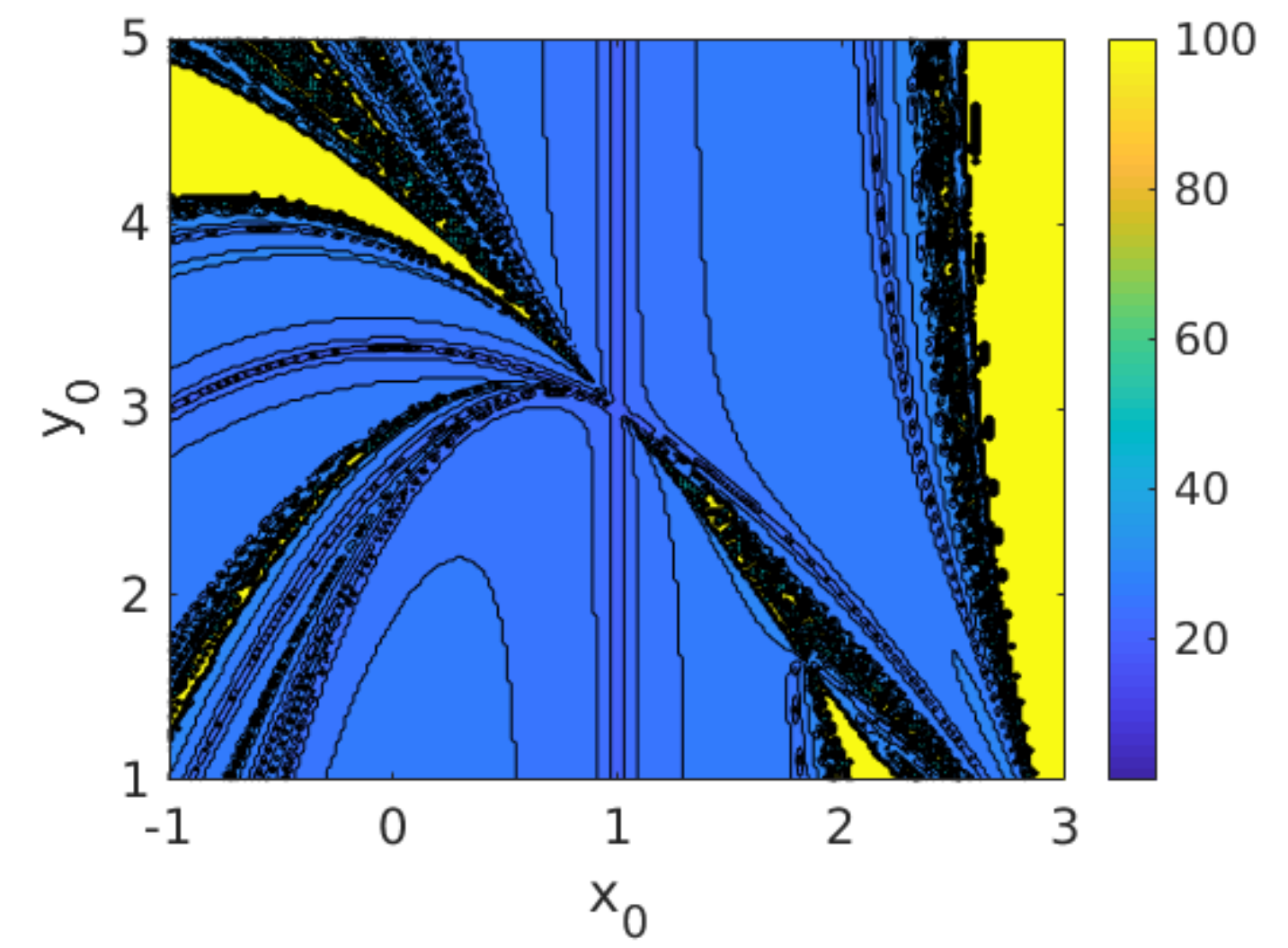}
\includegraphics[trim = 90pt 240pt 110pt 250pt,clip = true, width=0.32\textwidth]
{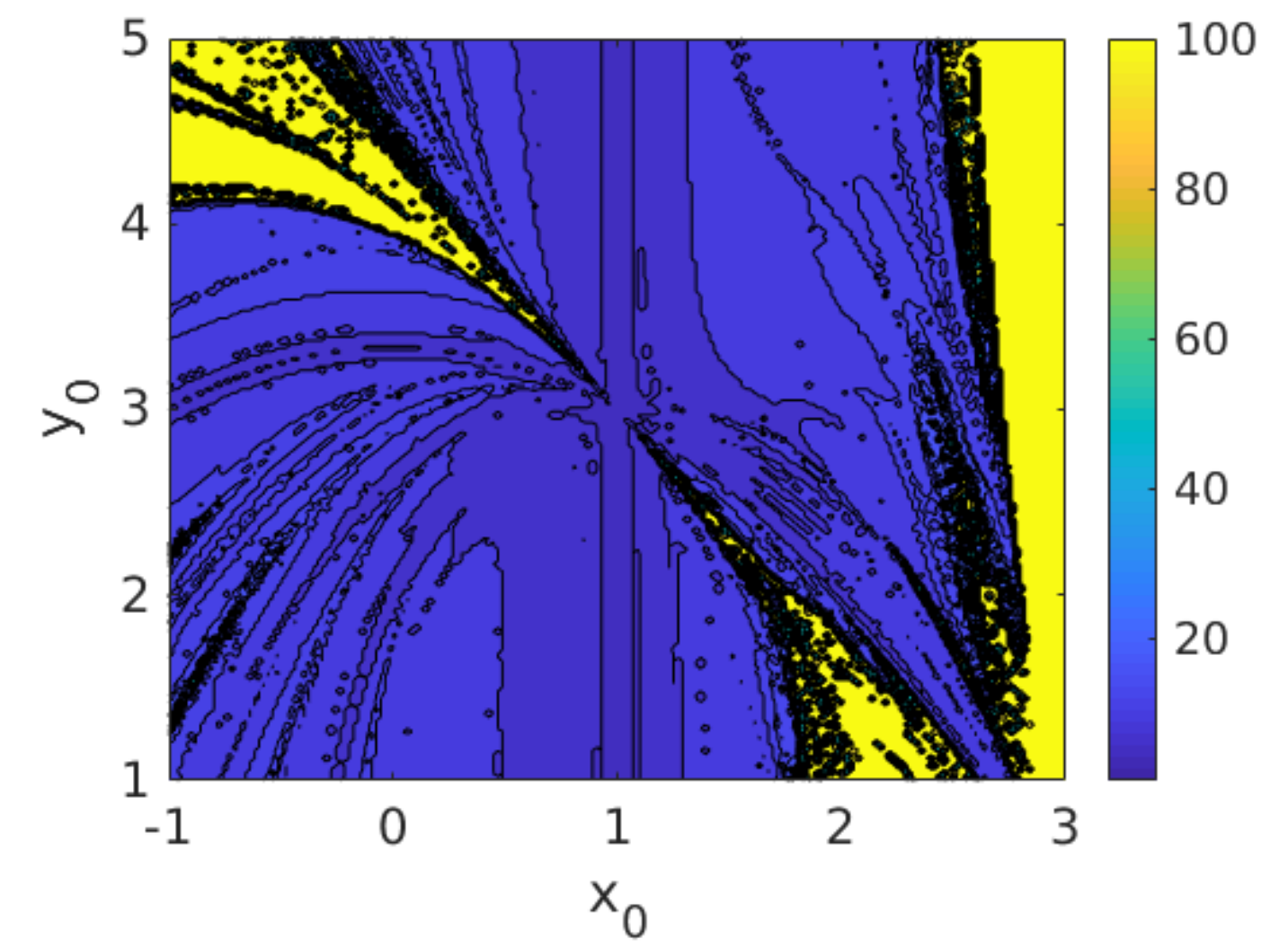}
\includegraphics[trim = 90pt 240pt 110pt 250pt,clip = true, width=0.32\textwidth]
{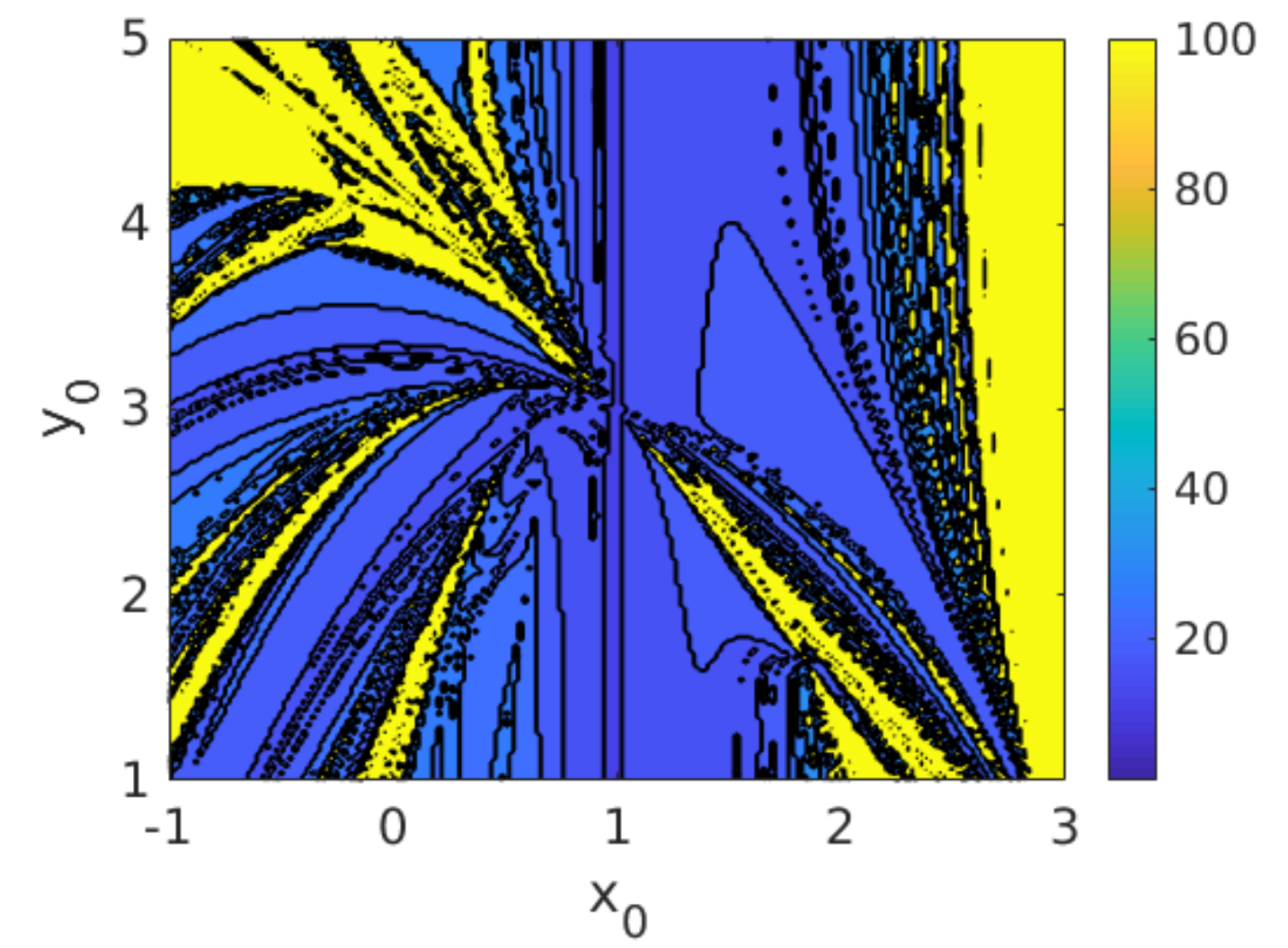}
\caption{Domain of convergence to $x_+=x_-$ for \eqref{eqn:poly2indef} with $\eps = 0$. 
Left: Newton iterations; 
center: Newton-Anderson(1); right: accelerated-Newton with $C = 0.5$ and 
$\alpha = 0.35$.
\label{fig:dofc-eps0}
}
\end{figure}
In this section, the domains of convergence for Newton, Newton-Andreson(1),
and accelerated-Newton are compared on standard benchmark
degenerate and nearly degenerate polynomial systems.

\begin{figure}[ht!]
\centering
\includegraphics[trim = 90pt 240pt 110pt 250pt,clip = true, width=0.32\textwidth]
{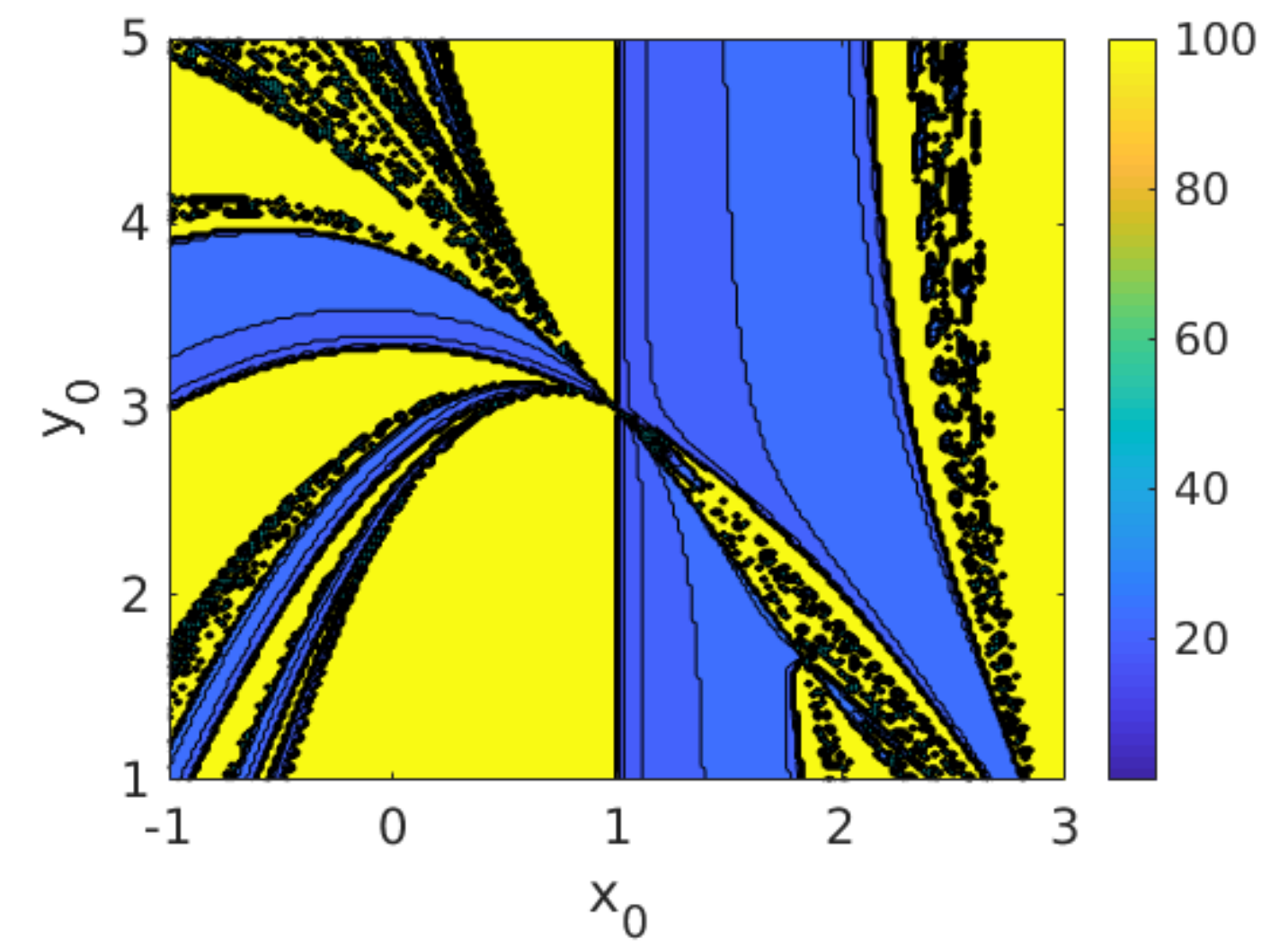}
\includegraphics[trim = 90pt 240pt 110pt 250pt,clip = true, width=0.32\textwidth]
{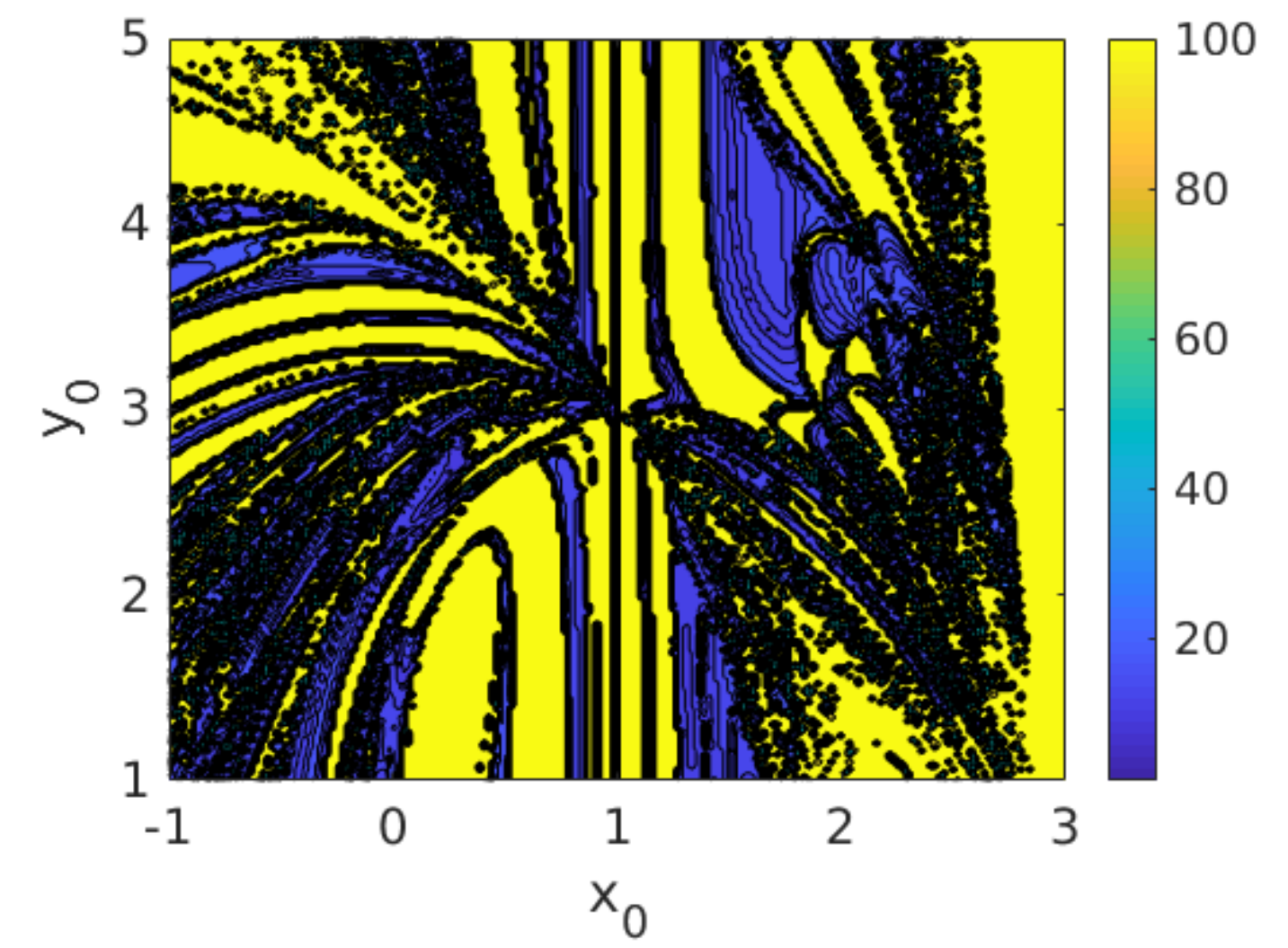}
\includegraphics[trim = 90pt 240pt 110pt 250pt,clip = true, width=0.32\textwidth]
{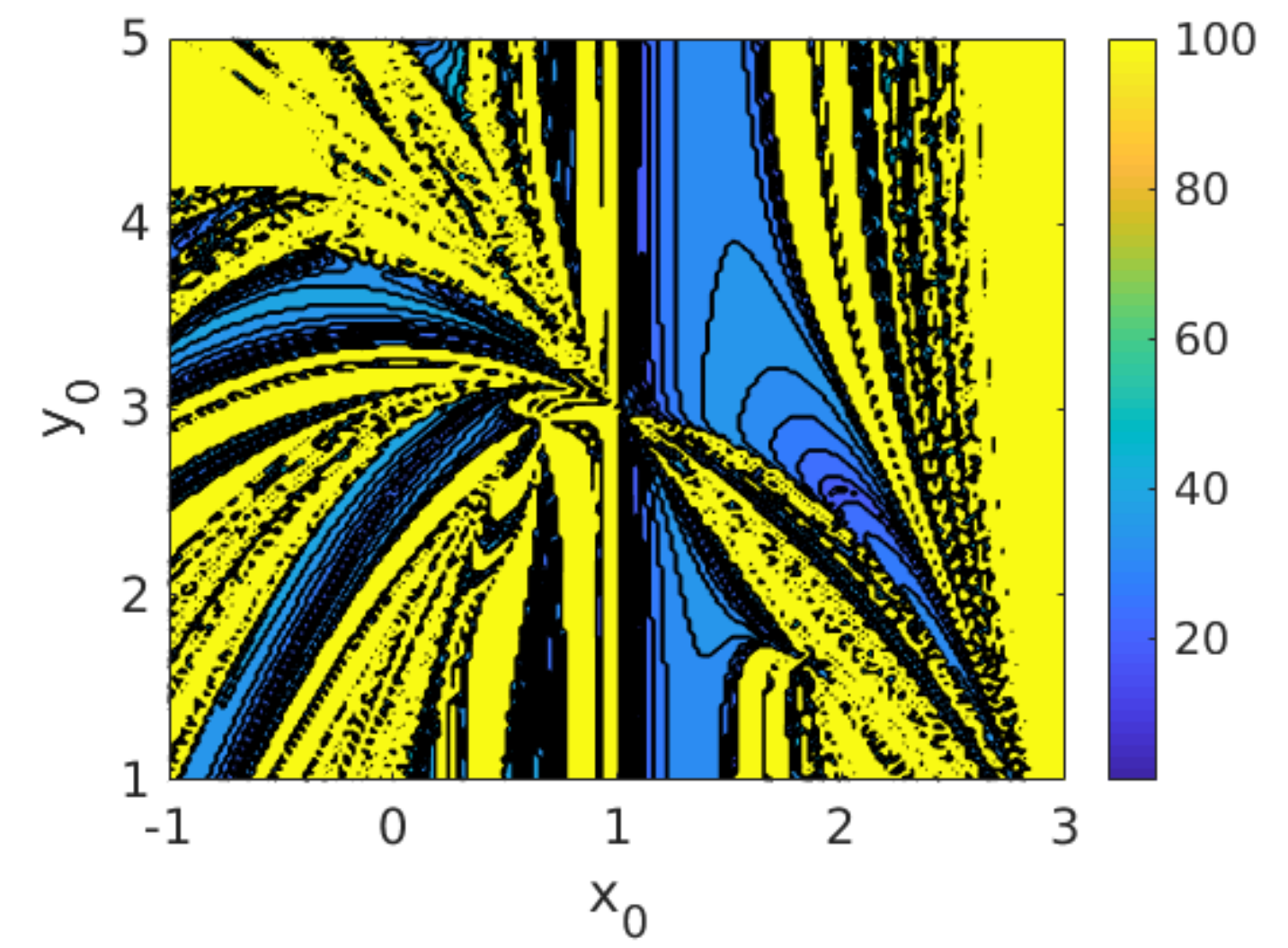}
\caption{Domain of convergence to $x_+$ for \eqref{eqn:poly2indef} 
with $\eps = 10^{-6}$. Left: Newton iterations; 
center: Newton-Anderson(1); right: accelerated-Newton.
\label{fig:dofc-eps-xp}
}
\end{figure}
\begin{figure}[ht!]
\centering
\includegraphics[trim = 90pt 240pt 110pt 250pt,clip = true, width=0.32\textwidth]
{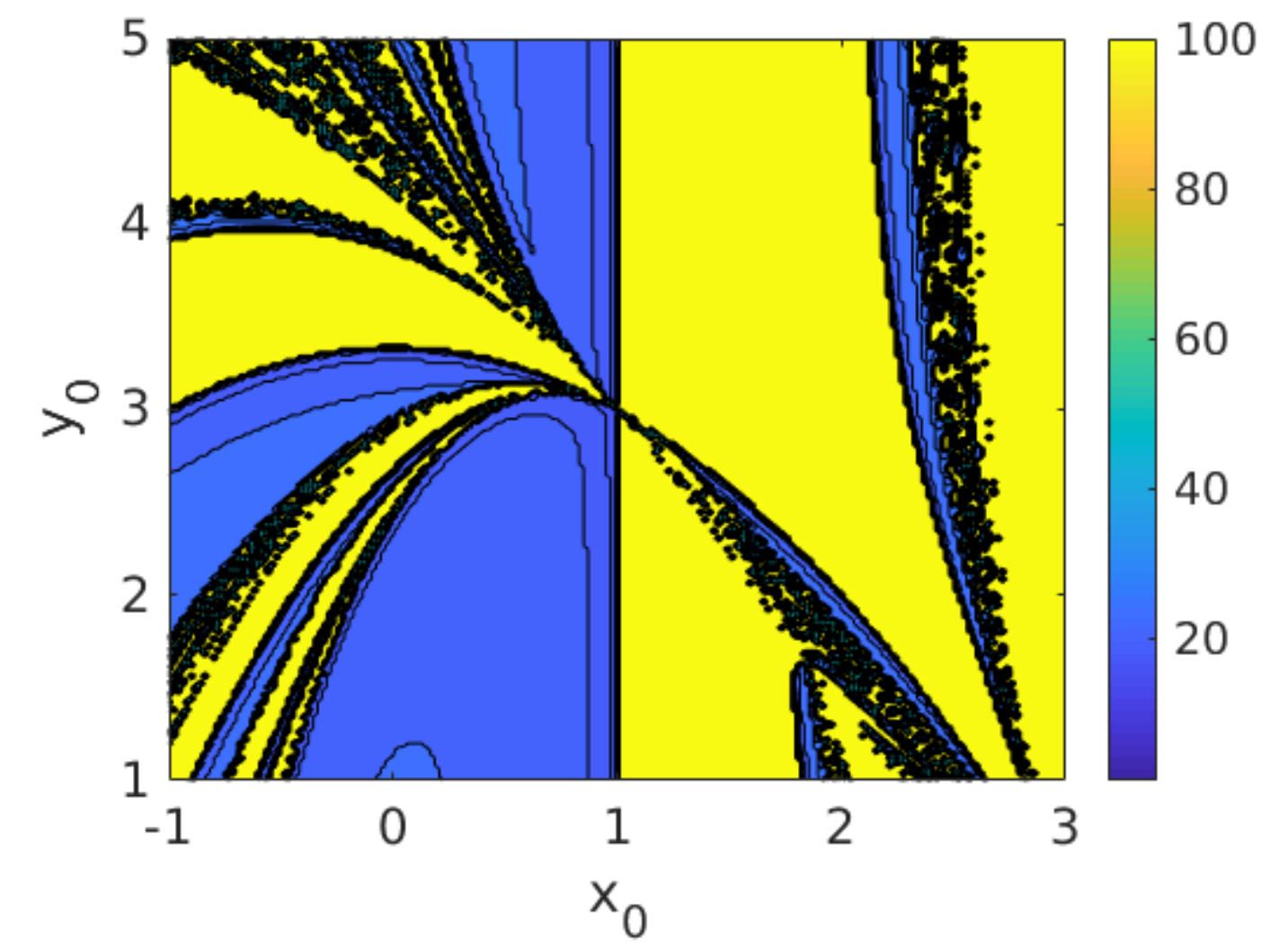}
\includegraphics[trim = 90pt 240pt 110pt 250pt,clip = true, width=0.32\textwidth]
{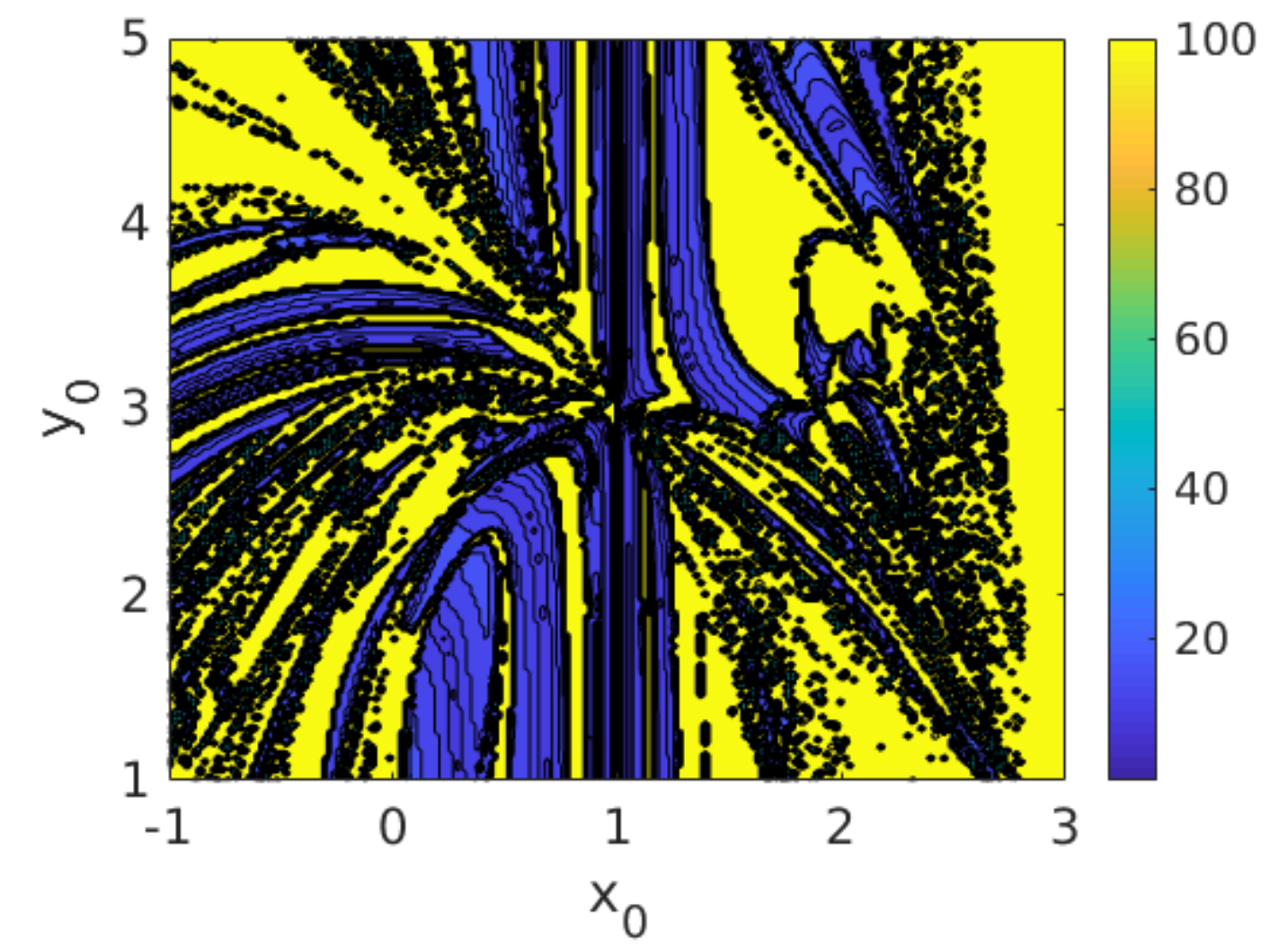}
\includegraphics[trim = 90pt 240pt 110pt 250pt,clip = true, width=0.32\textwidth]
{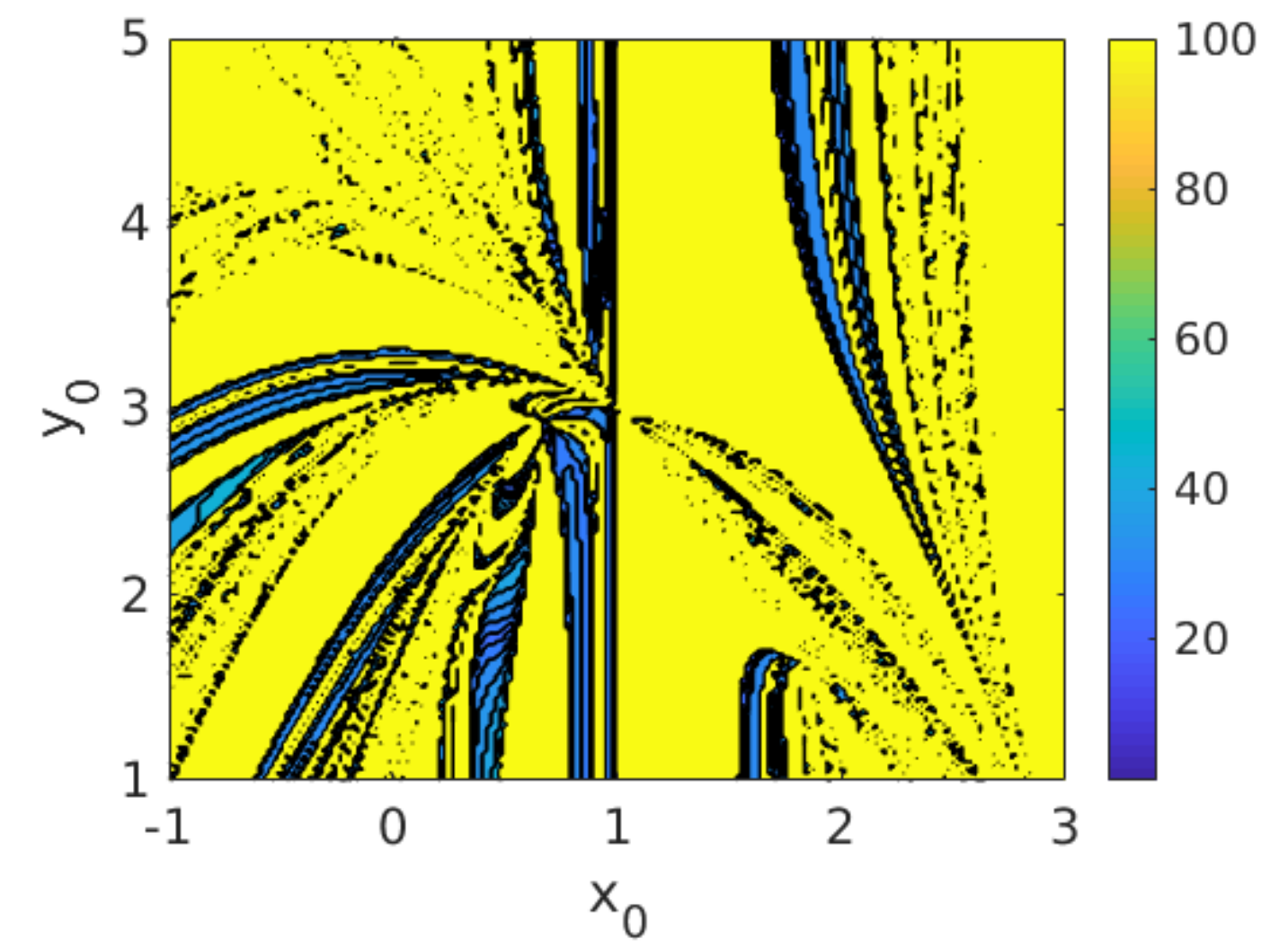}
\caption{Domain of convergence to $x_-$ for \eqref{eqn:poly2indef} 
with $\eps = 10^{-6}$. 
Left: Newton iterations; center: Newton-Anderson(1); right: accelerated-Newton.
\label{fig:dofc-eps-xm}
}
\end{figure}
The first investigated problem is the polynomial system $f: \R^2 \goto \R^2$ from 
\cite[Section 5]{DeKe85}.
\begin{align}\label{eqn:poly2indef}
f(x_1,x_2) = \begin{pmatrix} (x_1-1) + (x_2-3)^2\\
\eps(x_2-3) + (3/2) (x_1-1)(x_2-3) + (x_2-3)^2 + (x_2-3)^3 \end{pmatrix}.
\end{align}

As presented in \cite{DeKe85}, for $\eps > 0$ there are two distinct roots
in the vicinity of $x_+ = (1,3)^T$, the second being $x_- = (1-\eta^2, 3+\eta)^T$,
for $\eta = 1 - \sqrt{1 + 2\eps}$. For $\eps=0$, the degenerate case,
the two roots coincide; and, for $\eps \ll 1$ but $\eps > 0$, the near-degenerate case, 
the two roots are distinct but close. 
For figures \ref{fig:dofc-eps0}-\ref{fig:dofc-eps-xm}, each algorithm is run from
the the intial iterate $(x_0,y_0)$ for $x_0 \in [-1,3]$ and $y_0 \in [1,5]$ over a 
grid of 40,000 equally spaced points, and the number of iterations to convergence to
a particular solution is displayed.
The maximum of 100 iterations indicates not converging to the specified solution.
As demonstrated in figure \ref{fig:dofc-eps0},
Algorithm \ref{alg:ander} has a similar domain of convergence to 
Newton (Algorithm \ref{alg:newton}) and accelerated-Newton, (Algorithm \ref{alg:KS}); 
but, as illustrated in figures \ref{fig:dofc-eps-xp}-\ref{fig:dofc-eps-xm}, 
as $\eps$ is varied away from zero, the Anderson method has substantially 
different domains of attraction to the roots $x_-$ and $x_+$. For this example,
Algorithm \ref{alg:KS} was run with parameters $C = 0.5$ and $\alpha = 0.35$,
presented as optimal for $\eps = 10^{-6}$ in \cite{DeKe85}.

The second system investigated for domain of convergence is 
\eqref{eqn:poly3indef}, from subsection \ref{subsec:degen}.  
In this system, the Jacobian degenerates
for coordinates $x_3 = 0$, and also for $x_1=x_2 = 0$.  Here the domain of convergence is 
investigated for $x_0 = (a, b, 10^{-3})$, as the first two components
are varied, again over a grid of 40,000 equally spaced initial iterates.  
This example points out that the domains of convergence are similar for 
the three algorithms,
but the accelerated-Newton method is sensitive to parameter choice.  
With $\alpha$ and $C$ chosen small enough, the convergence is similar to Newton and
Newton-Anderson, but with the parameters chosen too large ($C=1, \alpha = 0.9$),
the star-shaped domain of convergence is considerably smaller.

\begin{figure}[ht!]
\label{fig:dofc-p3}
\centering
\includegraphics[trim = 90pt 240pt 110pt 250pt,clip = true, width=0.35\textwidth]
{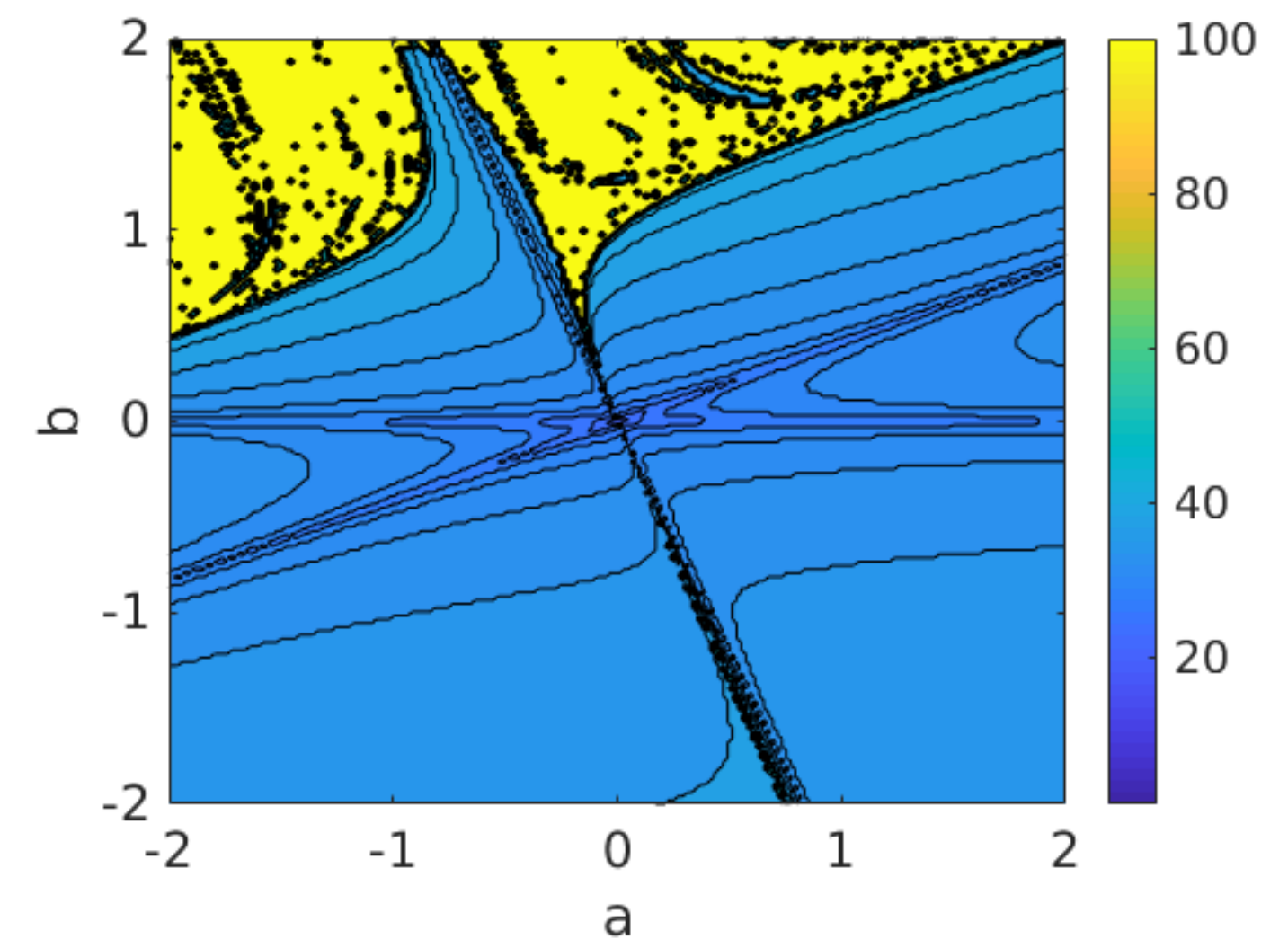}~\hfil~
\includegraphics[trim = 90pt 240pt 110pt 250pt,clip = true, width=0.35\textwidth]
{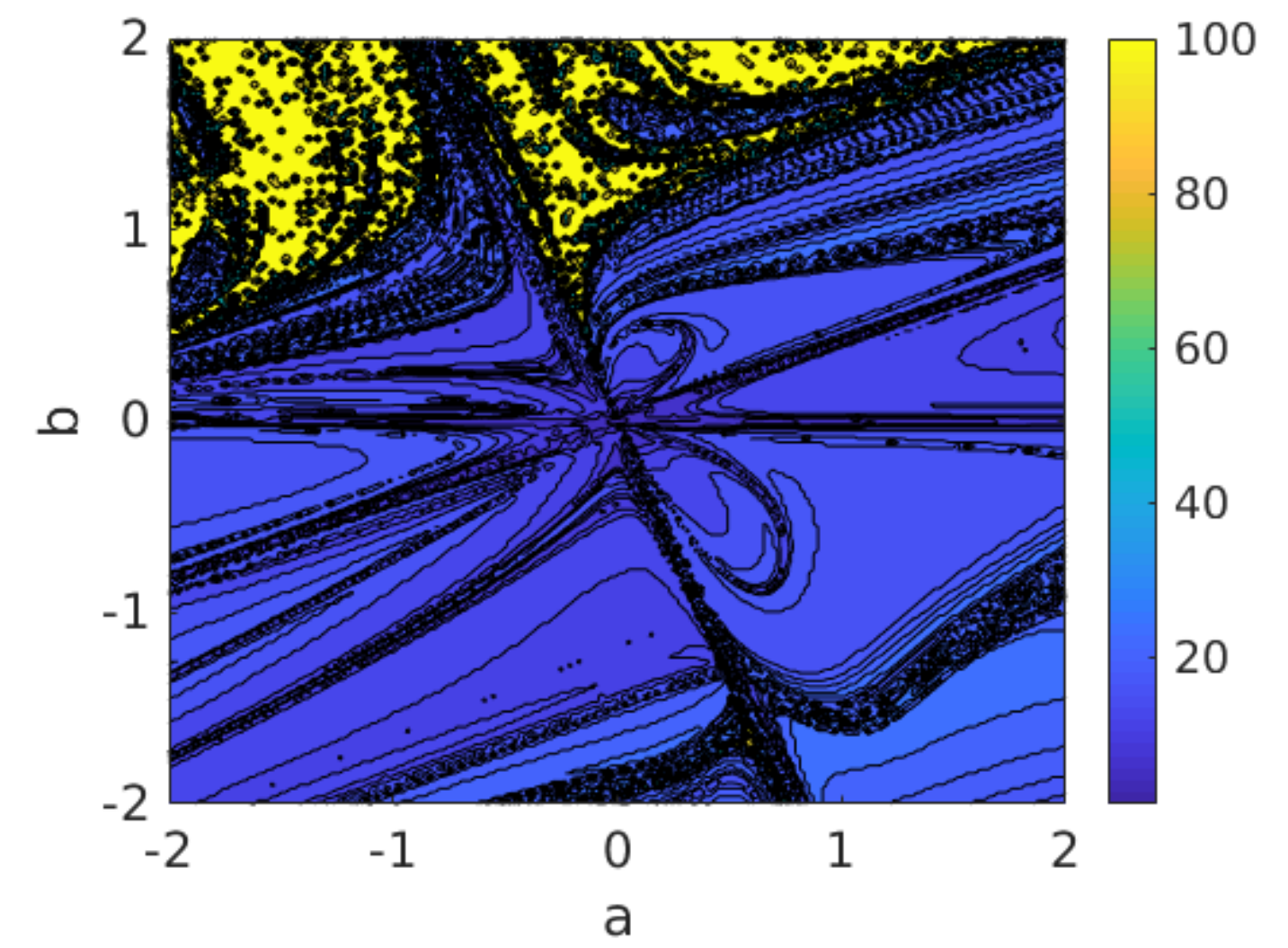}\\
\includegraphics[trim = 90pt 240pt 110pt 250pt,clip = true, width=0.35\textwidth]
{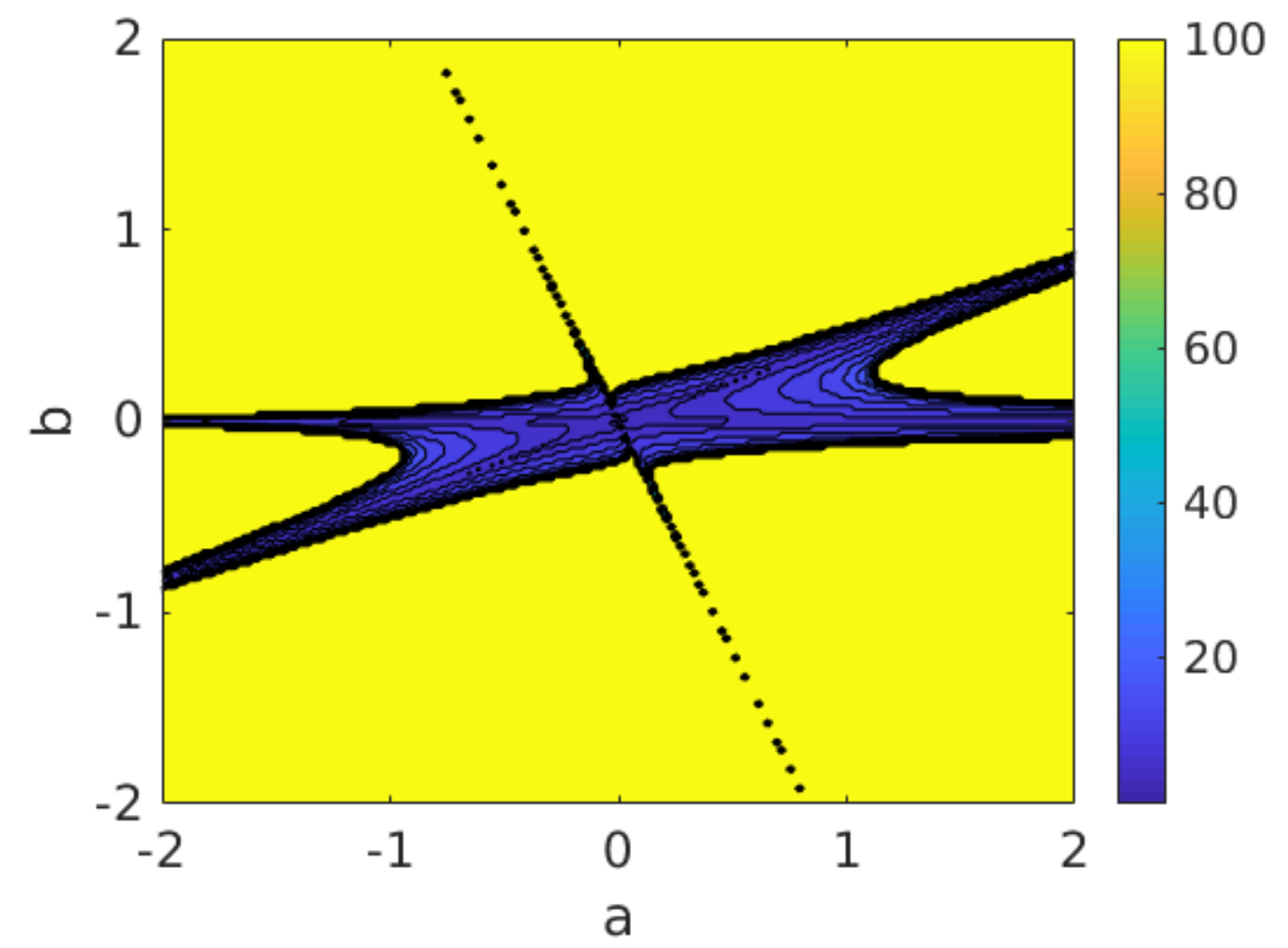}~\hfil~
\includegraphics[trim = 90pt 240pt 110pt 250pt,clip = true, width=0.35\textwidth]
{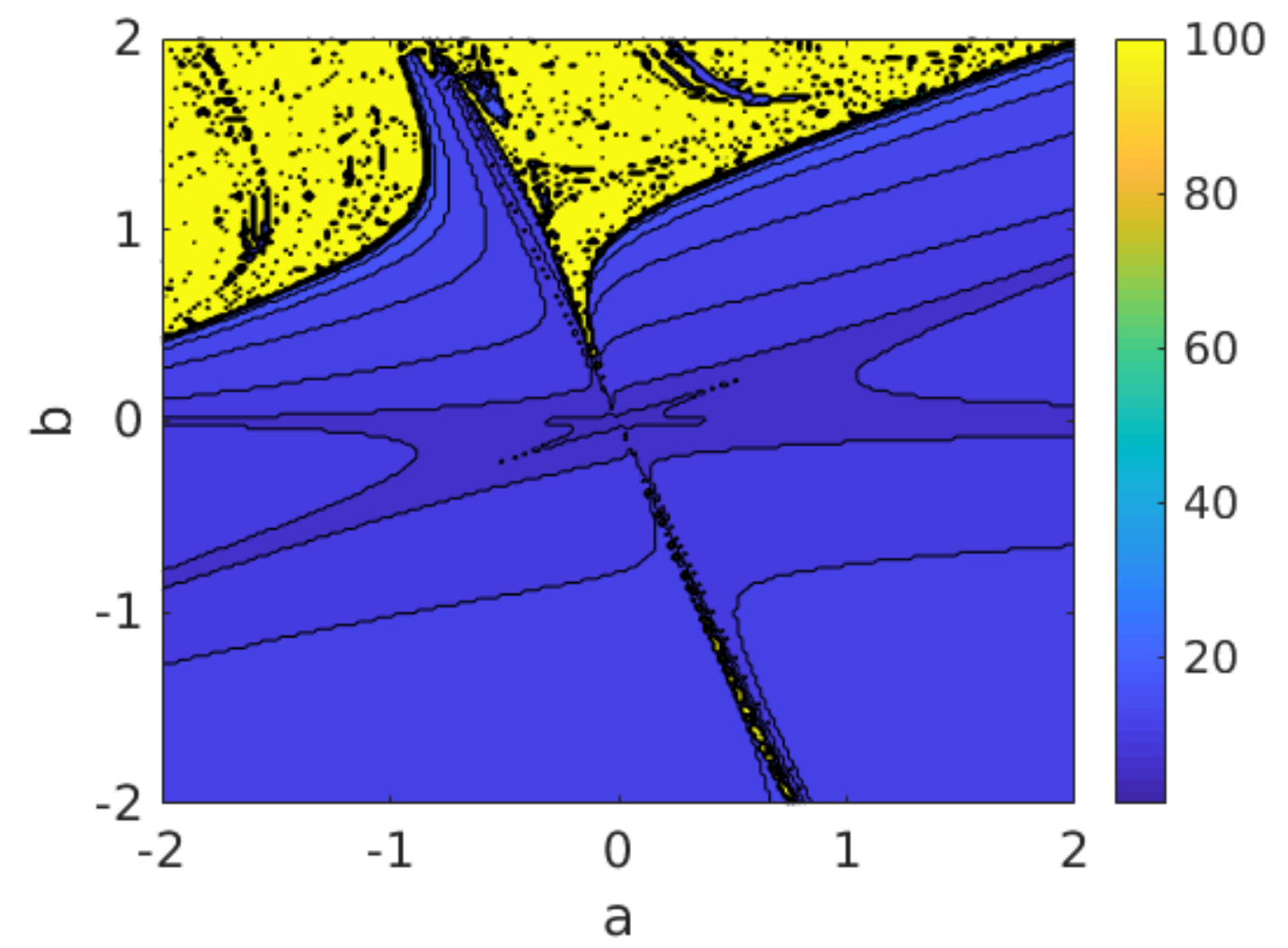}
\caption{Domain of convergence to $x^\ast = (0,0,0)^T$ 
for \eqref{eqn:poly3indef}, with
$x_0 = (a,b,10^{-3})^T$. Top left: Newton iterations; 
top right: Newton-Anderson(1); 
bottom left: accelerated-Newton with $C = 1.0, \alpha = 0.9$; 
bottom right: accelerated-Newton with $C = 0.1, \alpha = 0.35$.}
\end{figure}

\section{Conclusion}
\label{sec:conc}
This results of this numerical benchmarking investigation of the Newton-Anderson
method indicate a superlinear but (usually) subquadratic order of convergence
on both degenerate problems and nondegenerate problems.
The superlinear convergence for nondegenerate problems is also shown theoretically.
The Newton-Anderson method is seen to be 
advantageous in both the degenerate setting and for problems in which Newton's
method has an extended linearly converging preasymptotic regime, such as the 
trigonometric function \eqref{eqn:trig} considered here.
Further investigations show the method can show a similar domain of convergence
to the standard Newton iteration for degenerate problems; 
however, under small perturbation of a problem parameter,
the Newton-Anderson method may be attracted to a different solution than the Newton 
iteration.

\vfil

\section*{Acknowledgements} SP was partially supported by NSF DMS 1852876.
The authors would like to thank David Burrell (University of Florida)
for preliminary numerical results on degenerate systems.

  \bibliographystyle{elsarticle-num} 
  \bibliography{ander}


\end{document}